
\def\LaTeX{%
  \let\Begin\begin
  \let\End\end
  \let\salta\relax
  \let\finqui\relax
  \let\futuro\relax}

\def\UK{\def\our{our}\let\sz s}
\def\USA{\def\our{or}\let\sz z}

\UK



\LaTeX

\USA


\salta

\documentclass[twoside,12pt]{article}
\setlength{\textheight}{24cm}
\setlength{\textwidth}{16cm}
\setlength{\oddsidemargin}{2mm}
\setlength{\evensidemargin}{2mm}
\setlength{\topmargin}{-15mm}
\parskip2mm


\usepackage{color}
\usepackage{amsmath}
\usepackage{amsthm}
\usepackage{amssymb}
\usepackage[mathcal]{euscript}

%
%


\definecolor{viola}{rgb}{0.3,0,0.7}
\definecolor{ciclamino}{rgb}{0.5,0,0.5}


\bibliographystyle{plain}


%
\newtheorem{theorem}{Theorem}[section]
\newtheorem{remark}[theorem]{Remark}

\newtheorem{definition}[theorem]{Definition}

\finqui

\def\Beq{\Begin{equation}}
\def\Eeq{\End{equation}}
\def\Bsist{\Begin{eqnarray}}
\def\Esist{\End{eqnarray}}

\def\Bthm{\Begin{theorem}}
\def\Ethm{\End{theorem}}
\def\Blem{\Begin{lemma}}
\def\Elem{\End{lemma}}

\def\Bdim{\Begin{proof}}
\def\Edim{\End{proof}}
\def\Bcenter{\Begin{center}}
\def\Ecenter{\End{center}}
\let\non\nonumber




\def\step #1 \par{\medskip\noindent{\bf #1.}\quad}


\def\Frechet{Fr\'echet}
\def\aand{\quad\hbox{and}\quad}


\def\generaliz{generali\sz}


\def\multibold #1{\def\arg{#1}%
  \ifx\arg\pto \let\next\relax
  \else
  \def\next{\expandafter
    \def\csname #1#1#1\endcsname{{\bf #1}}%
    \multibold}%
  \fi \next}

\def\pto{.}

\def\multical #1{\def\arg{#1}%
  \ifx\arg\pto \let\next\relax
  \else
  \def\next{\expandafter
    \def\csname cal#1\endcsname{{\cal #1}}%
    \multical}%
  \fi \next}


\def\multimathop #1 {\def\arg{#1}%
  \ifx\arg\pto \let\next\relax
  \else
  \def\next{\expandafter
    \def\csname #1\endcsname{\mathop{\rm #1}\nolimits}%
    \multimathop}%
  \fi \next}

\multibold
qwertyuiopasdfghjklzxcvbnmQWERTYUIOPASDFGHJKLZXCVBNM.

\multical
QWERTYUIOPASDFGHJKLZXCVBNM.

\multimathop
dist div dom meas sign supp .


\def\Accorpa #1#2 #3 {\gdef #1{\eqref{#2}--\eqref{#3}}%
  \wlog{}\wlog{\string #1 -> #2 - #3}\wlog{}}


\def\infess{\mathop{\rm inf\,ess}}
\def\supess{\mathop{\rm sup\,ess}}

\def\neto{\mathrel{{\scriptscriptstyle\nearrow}}}
\def\seto{\mathrel{{\scriptscriptstyle\searrow}}}

\def\graffe #1{\mathopen\{#1\mathclose\}}

\def\<#1>{\mathopen\langle #1\mathclose\rangle}
\def\norma #1{\mathopen \| #1\mathclose \|}

\def\iO{\int_\Omega}
\def\iG{\int_\Gamma}

\def\dt{\partial_t}
\def\dn{\partial_{\bf n}}

\def\checkmmode #1{\relax\ifmmode\hbox{#1}\else{#1}\fi}

\def\aet{\checkmmode{a.e.\ in~$(0,T)$}}

\def\aat{\checkmmode{for a.a.~$t\in(0,T)$}}


\def\erre{{\mathbb{R}}}

\def\enne{{\mathbb{N}}}




\def\genspazio #1#2#3#4#5{#1^{#2}(#5,#4;#3)}
\def\spazio #1#2#3{\genspazio {#1}{#2}{#3}T0}

\def\L {\spazio L}
\def\H {\spazio H}
\def\W {\spazio W}


\def\Hx #1{H^{#1}(\Omega)}

\def\HxG #1{H^{#1}(\Gamma)}

\def\Cx #1{C^{#1}(\overline\Omega)}

\def\LS #1{L^{#1}(\Sigma)}

\def\Hdue{\Hx 2}

\def\HdueG{\HxG 2}



\let\theta\vartheta

\let\phi\varphi

\let\TeXchi\chi                         
\newbox\chibox
\setbox0 \hbox{\mathsurround0pt $\TeXchi$}
\setbox\chibox \hbox{\raise\dp0 \box 0 }
\def\chi{\copy\chibox}



\def\suG{{\vrule height 5pt depth 4pt\,}_\Gamma}

\def\fG{f_\Gamma}
\def\yG{y_\Gamma}
\def\uG{u_\Gamma}
\def\vG{v_\Gamma}
\def\hG{h_\Gamma}
\def\xiG{\xi_\Gamma}
\def\qG{q_\Gamma}
\def\uGmin{u_{\Gamma,{\rm min}}}
\def\uGmax{u_{\Gamma,{\rm max}}}

\def\yz{y_0}

\def\Mz{M_0}
\def\ustar{u_*}
\def\vstar{v_*}

\def\bQ{b_Q}
\def\bS{b_\Sigma}
\def\bO{b_\Omega}
\def\bG{b_\Gamma}
\def\bz{b_0}

\def\rmin{r_-}
\def\rmax{r_+}

\def\zQ{z_Q}
\def\zS{z_\Sigma}
\def\zO{z_\Omega}
\def\zG{z_\Gamma}

\def\Uad{\calU_{ad}}
\def\uopt{\bar u_\Gamma}
\def\yopt{\bar y}
\def\yGopt{\bar y_\Gamma}

\def\redJ{\widetilde\calJ}

\def\vO{v^\Omega}

\def\VG{V_\Gamma}
\def\HG{H_\Gamma}
\def\Vp{V^*}

\def\normaV #1{\norma{#1}_V}
\def\normaH #1{\norma{#1}_H}

\def\normaVp #1{\norma{#1}_*}

\def\nablaG{\nabla_{\!\Gamma}}
\def\DeltaG{\Delta_\Gamma}

\let\hat\widehat

\def\Pi{\hat\pi}

\def\mz{m_0}


\Begin{document}


\title{Second-order analysis of a boundary control\\[0.3cm] 
problem for the viscous Cahn--Hilliard equation\\[0.3cm]
  with dynamic boundary condition
}

\author{}
\date{}
\maketitle
\Bcenter
\vskip-1cm
{\large\sc Pierluigi Colli$^{(1)}$}\\
{\normalsize e-mail: {\tt pierluigi.colli@unipv.it}}\\[.25cm]
{\large\sc M. Hassan Farshbaf-Shaker$^{(2)}$}\\
{\normalsize e-mail: {\tt Hassan.Farshbaf-Shaker@wias-berlin.de}}\\[.25cm]
{\large\sc Gianni Gilardi$^{(1)}$}\\
{\normalsize e-mail: {\tt gianni.gilardi@unipv.it}}\\[.25cm]
{\large\sc J\"urgen Sprekels$^{(2,3)}$}\\
{\normalsize e-mail: {\tt sprekels@wias-berlin.de}}\\[.45cm]
$^{(1)}$
{\small Dipartimento di Matematica ``F. Casorati'', Universit\`a di Pavia}\\
{\small via Ferrata 1, 27100 Pavia, Italy}\\[.2cm]
$^{(2)}$
{\small Weierstrass Institute for Applied Analysis and Stochastics}\\
{\small Mohrenstrasse 39, 10117 Berlin, Germany}\\[2mm]
$^{(3)}$
{\small Department of Mathematics}\\
{\small Humboldt-Universit\"at zu Berlin}\\
{\small Unter den Linden 6, 10099 Berlin, Germany}\\
[1cm]
\Ecenter
\begin{center}
{\em Dedicated to the memory of Prof.\ Dr.\ Viorel Arn\u{a}utu}
\end{center}
\Begin{abstract}
In this paper we establish second-order sufficient optimality conditions 
for a boundary control problem 
that has been introduced and studied by three of the authors in the
preprint arXiv:1407.3916.
This control problem regards
the viscous Cahn--Hilliard equation 
with possibly singular potentials and dynamic boundary conditions.
\vskip3mm

\noindent {\bf Key words:}
Cahn--Hilliard equation, dynamic boundary conditions, phase separation,
singular potentials, optimal control, first and second order optimality conditions,
adjoint state system.
\vskip3mm
\noindent {\bf AMS (MOS) Subject Classification:} 35K55 (35K50, 82C26)
\End{abstract}

\salta

\pagestyle{myheadings}
\newcommand\testopari{\sc Colli \ --- \ Farshbaf-Shaker \ --- \ Gilardi \ --- \ Sprekels}
\newcommand\testodispari{\sc Second-order analysis of control for the viscous C.--H. equation}
\markboth{\testodispari}{\testopari}

\finqui


\section{Introduction}
\label{Intro}
\setcounter{equation}{0}
This paper deals with second-order optimality conditions of a special boundary control problem for the viscous Cahn--Hilliard equation with dynamic boundary conditions. It continues the work \cite{CGSopt} by three of the present authors in which the first-order necessary conditions of optimality were derived. For the work of other authors concerning the optimal
control of Cahn--Hilliard systems, we refer the reader to the references given in~\cite{CGSopt}.

Crucial contributions in \cite{CGSopt} were the derivation of the adjoint problem, whose form turned out to be nonstandard, and an existence result for its solutions. As is well known, first-order conditions are in the case of nonlinear equations usually not sufficient for optimality. Also, second-order sufficient optimality conditions for nonlinear optimal control problems are essential both in the numerical analysis and for the construction of reliable optimization algorithms. For instance, the strong convergence of optimal controls and states for numerical discretizations of the problem rests heavily on the availability of second-order sufficient optimality conditions; furthermore, one can show that numerical algorithms such as SQP methods are locally convergent if second-order sufficient optimality conditions hold true. For a general discussion of second-order sufficient conditions for elliptic and parabolic control problems we refer to \cite{Tr} and references therein; for the case of control problems involving phase field models we refer to, e.\,g., \cite{CS,Heik}. 

In this paper, we aim to
establish second-order sufficient optimality conditions for the boundary control problem studied in \cite{CGSopt}.   
To this end, we assume that an open, bounded and connected set $\Omega\subset\erre^3$, with smooth boundary $\Gamma$
and unit outward normal~{\bf n}, and some final time $T>0$ are given, and we set $Q:=\Omega\times (0,T)$ and
$\Sigma:=\Gamma\times (0,T)$. Moreover, we denote by $\Delta_\Gamma$, $\nabla_\Gamma$, $\partial_{\bf n}$, the
Laplace--Beltrami operator, the surface gradient, and the outward normal derivative on $\Gamma$, in this order. 
We make the following general assumptions:

\vspace{3mm}\noindent
{\bf (A1)} \quad\,There are given nonnegative constants $\bQ,\bS,\bO,\bG,\bz$, which do not all vanish, functions 
$z_Q\in L^2(Q)$, $z_\Sigma\in L^2(\Sigma)$, $z_\Omega\in L^2(\Omega)$, $z_\Gamma\in L^2(\Gamma)$, as well as a constant $\,M_0>0\,$ and functions $\,\uGmin\in\LS\infty\,$ and $\,\uGmax\in\LS\infty\,$ with $\,\uGmin\leq\uGmax$ a.\,e. in 
$\Sigma$.

\vspace{3mm}\noindent
{\bf (A2)} \quad\,There are given constants $\,-\infty\le r_-<0<r_+\le +\infty\,$ and two functions 
$\,f,f_\Gamma:(r_-,r_+)\to [0,+\infty)$\, such that the following holds:
\begin{align}
\label{conf1}
&f,f_\Gamma\in C^4(r_-,r_+), \quad f(0)=f_\Gamma(0)=0,\\[1mm]
\label{conf2}
&f'' \,\mbox{ and }\, f_\Gamma'' \,\mbox{ are bounded from below},\\[1mm]
\label{conf3}
& \lim\limits_{r\seto\rmin} f'(r)
  = \lim\limits_{r\seto\rmin} f_\Gamma'(r) 
  = -\infty 
  \aand
  \lim\limits_{r\neto\rmax} f'(r)
  = \lim\limits_{r\neto\rmax} \fG'(r) 
  = +\infty\,,\\[1mm]
\label{conf4}
& |f'(r)| \leq \eta \,|\fG'(r)| + C
  \quad \hbox{for some $\eta,\, C>0$ and every $r\in(\rmin,\rmax)$}.
\end{align} 

\noindent
In fact, \eqref{conf1} is fully used only in the last part of the paper
and many of our results hold under a weaker assumption.
We also note that the conditions (\ref{conf1})--(\ref{conf4}) allow for the possibility
of splitting $f'$ in (\ref{conf3}) in the form $f'=\beta+\pi$,
where $\beta$ is a monotone function that diverges at~$r_\pm$
and $\pi$ is a perturbation with a bounded derivative. Since the same is true for~$\fG$,
the general  assumptions of \cite{CGS} are satisfied.
Typical and important examples for $f$ and $f_\Gamma$ are 
the classical regular potential $f_{\rm reg}$ 
and the logarithmic double-well potential $f_{\rm log}$ given~by
\begin{align}
  & f_{\rm reg}(r) = \frac14(r^2-1)^2 \,,
  \quad r \in \erre 
  \label{regpot}
  \\
  & f_{\rm log}(r) = ((1+r)\ln (1+r)+(1-r)\ln (1-r)) - c r^2 \,,
  \quad r \in (-1,1),
  \label{logpot}
\end{align}
where in the latter case we assume that $c>0$ is so large that $f_{\rm log}$ is nonconvex.

With the above assumptions, we consider the following tracking type optimal boundary
control problem:

\noindent ({\bf CP}) Minimize  
\Bsist
  & \calJ(y,\yG,\uG)
  & := \frac\bQ 2 \, \norma{y-\zQ}_{L^2(Q)}^2
  + \frac\bS 2 \, \norma{y_\Gamma-\zS}_{L^2(\Sigma)}^2
  + \frac\bO 2 \, \norma{y(T)-\zO}_{L^2(\Omega)}^2
  \non
  \\
  && \quad +\, \frac\bG 2 \, \norma{\yG(T)-\zG}_{L^2(\Gamma)}^2
  + \frac\bz 2 \, \norma\uG_{L^2(\Sigma)}^2
  \label{Icost}
\Esist
subject to the control constraint  
\Bsist
  & \uG\in\Uad := 
  & \bigl\{ \vG\in H^1(0,T;L^2(\Gamma))\cap\LS\infty:
  \non
  \\
  && \,\,\,\,\uGmin\leq\vG\leq\uGmax \,\,\,\mbox{a.\,e. on }\,\Sigma,\,\, \norma{\dt\vG}_2\leq\Mz
  \bigr\}  
  \label{Iuad}
\Esist
and to the Cahn--Hilliard equation with nonlinear dynamic boundary conditions as the state system, 
\Bsist
  & \dt y - \Delta w = 0   \quad \mbox{in }\, Q,
  \label{Iprima}
  \\[0.2cm]
  & w =  \dt y - \Delta y + f'(y)
  \quad \mbox{in }\,Q,
  \label{Iseconda}
  \\[0.2cm]
  & \partial_{\bf n} w = 0
  \quad \mbox{on }\, \Sigma,
  \label{Ibc}
  \\[0.2cm]
  & \yG = y\suG\quad \mbox{on }\, \Sigma,
  \label{ICorrespond}
  \\[0.2cm]
  &\dt\yG + \dn y - \DeltaG\yG + f'_\Gamma(\yG) = \uG
  \quad \hbox{on $\, \Sigma$,}
  \label{Iterza}
  \\[0.2cm]
  & y(\cdot,0) = \yz \quad \hbox{in $\, \Omega$},\quad \yG(\cdot,0) = y_{0_\Gamma} \quad \hbox{on $\, \Gamma$.}
  \label{Icauchy}
\Esist

\noindent Here, and throughout this paper, we generally assume that the admissible set $\Uad$ is nonempty. Moreover, we
postulate:

\vspace{3mm}
\noindent
{\bf (A3)} \quad\,$y_0\in H^2(\Omega)$, $y_{0_\Gamma}:=y_{0|\Gamma}\in H^2(\Gamma)$, and it holds
(notice that $\yz\in\Cx0$)
\begin{equation}
\label{boundy0}
r_- < y_0 < r_+ \quad\mbox{in $\overline\Omega$}.
\end{equation}

\vspace{2mm}\noindent
We remark at this place that in \cite{CGS} the additional assumption $\partial_{\bf n}
y_0=0$ was made; this postulate is however unnecessary for the results of \cite{CGS}
to hold, since it is nowhere used in the proofs.

\vspace{2mm}
The system (\ref{Iprima})--(\ref{Icauchy}) is an initial-boundary 
value problem with nonlinear dynamic boundary condition for a Cahn--Hilliard equation. 
In this connection, the unknown $y$ usually stands for the order parameter of an isothermal phase transition, 
and $w$ denotes the chemical potential of the system. 

Our paper is organized as follows: in Section \ref{STATEMENT}, we provide and collect some results 
proved in \cite{CGSopt,CGS} concerning the state system, 
and we study a certain linear counterpart thereof that will be employed repeatedly in the later analysis. 
In Section~3, the existence of the second-order Fr\'echet derivative of the control-to-state mapping will be shown. 
Section~4 then brings the derivation of the second-order sufficient condition of optimality.    

In order to simplify notation, we will in the following write $y_\Gamma$
for the trace $y_{|\Gamma}$ of a function $y\in H^1(\Omega)$
on $\Gamma$, and we introduce the abbreviations
\begin{align}
\label{spaces}
&V:=H^1(\Omega),\quad H:=L^2(\Omega), \quad V_\Gamma:=H^1(\Gamma), \quad
H_\Gamma:=L^2(\Gamma), \quad {\cal H}:=H\times H_\Gamma,\non\\[1mm]
&{\cal V}:=\{(v,v_\Gamma)\in V\times V_\Gamma:\,v_\Gamma=v_{|\Gamma}\},
\quad {\cal G}:=H^2(\Omega)\times H^2(\Gamma),\non\\[1mm]
& {\cal X}:=H^1(0,T;H_\Gamma)\cap
L^\infty(\Sigma), \quad {\cal Y}:= H^1(0,T;{\cal H})\cap L^\infty(0,T;{\cal V}),
\end{align}
and endow these spaces with their natural norms. Moreover, for the
 generic Banach space $X$ we denote by $X^*$ its dual space and by
$\|\cdot\|_X$ its norm. Furthermore, the symbol $\,\langle\,\cdot
\,,\,\cdot\,\rangle\,$ stands for the duality pairing between the spaces
$V^*$ and $V$, where it is understood that $H$ is embedded in $V^*$ in
the usual way, i.\,e., such that we have $\,\langle u,v\rangle=
(u,v)\,$ for every $u\in H$ and $v\in V$ with the standard inner product
$(\,\cdot\,,\,\cdot\,)$ of $H$. Finally, for $u\in V^*$ and $v\in
L^1(0,T;V^*)$  we define their generalized mean values $u^\Omega\in\erre$ 
and $v^\Omega \in L^1(0,T)$, respectively, by setting
\begin{equation}
\label{media}
u^\Omega :=\frac 1 {|\Omega|} \,\langle u,1\rangle \quad\,\mbox{and }\,
v^\Omega(t):=(v(t))^\Omega \quad\mbox{for a.\,e. }\,t\in (0,T),
\end{equation}
where $\,|\Omega|\,$ stands for the Lebesgue measure of $\Omega$.

 During the course of our analysis, we will make repeated use of the elementary Young's inequality
\Beq
  ab \leq \delta a^2 + \frac 1 {4\delta} \, b^2
  \quad \hbox{for every $a,b\geq 0$ and $\delta>0$,}
  \label{young}
\Eeq
of H\"older's inequality, and of Poincar\'e's inequality
\Beq
  \normaV v^2 \leq \widehat C \bigl( \normaH{\nabla v}^2 + |\vO|^2 \bigr)
  \quad \hbox{for every $v\in V$},
  \label{poincare}
\Eeq
where $\widehat C>0$~depends only on~$\Omega$.

Next, we recall a tool that is commonly used 
in the context of problems related to the Cahn--Hilliard equations.
We define
\Beq
  \dom\calN := \graffe{\vstar\in\Vp: \ \vstar^\Omega = 0}
  \aand
  \calN : \dom\calN \to \graffe{v \in V : \ \vO = 0}
  \label{predefN}
\Eeq
by setting, for $\vstar\in\dom\calN$,
\Beq
  {\calN\vstar \in V, \quad
  (\calN\vstar)^\Omega = 0 ,
  \aand
  \iO \nabla\calN\vstar \cdot \nabla z\,dx = \< \vstar , z >
  \quad \hbox{for every $z\in V$}}\,.
  \label{defN}
\Eeq
That is, $\calN\vstar$ is the unique solution $v$ to the \generaliz ed Neumann problem for $-\Delta$ 
with datum~$\vstar$ that satisfies~$\vO=0$.
Indeed, if $\vstar\in H$, then the above variational equation means that
$-\Delta {\cal N}\vstar = \vstar$ in $\Omega$ and $\dn\calN\vstar = 0$ on $\Gamma$.
Moreover, we have
\Beq
  \< \ustar , \calN \vstar >
  = \< \vstar , \calN \ustar >
  = \iO (\nabla\calN\ustar) \cdot (\nabla\calN\vstar)\,dx
  \quad \mbox{for all }\,\ustar,\vstar\in\dom\calN,
  \label{simmN}
\Eeq
whence also
\Beq
  2 \< \dt\vstar(t) , \calN\vstar(t) >
  = \frac d{dt} \iO |\nabla\calN\vstar(t)|^2\,dx
  = \frac d{dt} \, \normaVp{\vstar(t)}^2
  \quad \aat
  \label{dtcalN}
\Eeq
for every $\vstar\in\H1\Vp$ satisfying $(\vstar)^\Omega=0$ \aet.

\section{The state equation}\label{STATEMENT}
\setcounter{equation}{0}
\noindent  
At first, we specify our notion of solution to the state system
(\ref{Iprima})--(\ref{Icauchy}).

\begin{definition}\label{DefSolu}
{\rm Suppose that the general assumptions} {\bf (A1)}--{\bf (A3)} {\em are fulfilled,
and let $\uG\in {\cal X}$ be given. 
By a} solution {\em to (\ref{Iprima})--(\ref{Icauchy}) we mean a  triple 
$(y,\yG,w)$ that satisfies }
\begin{align}
  & y \in \W{1,\infty}H \cap \H1V \cap \L\infty\Hdue,
  \label{regy}
  \\[1mm]
  & \yG \in \W{1,\infty}\HG \cap \H1\VG \cap \L\infty\HdueG,
  \label{regyG}
  \\[1mm]
  & \yG(t) = y(t)\suG
  \quad \mbox{\rm for a.\,a. }\,t\in (0,T),
  \label{tracciay}
  \\[1mm]
  & 
  \rmin < \infess\limits_Q y \leq \supess\limits_Q y < \rmax\,,\quad
  \rmin < \infess\limits_\Sigma y_\Gamma \leq \supess\limits_\Sigma y_\Gamma
  < \rmax\,, 
  \label{tuttolinfty}
  \\[1mm]
  & w \in \L\infty\Hdue,
  \label{regw}
\end{align}
{\rm as well as, for almost every $t\in (0,T)$, the variational equations}
\Bsist
  && \iO \dt y(t) \, v \,dx
  + \iO \nabla w(t) \cdot \nabla v\,dx = 0,
  \label{prima}
  \\
  \noalign{\smallskip}
  && \iO w(t) \, v \,dx
  = \iO \dt y(t) \, v\,dx
  + \iG \dt\yG(t) \, \vG\,d\Gamma
  + \iO \nabla y(t) \cdot \nabla v\,dx
  \non\\[1mm]
  &&+ \iG \nablaG\yG(t) \cdot \nablaG\vG\,d\Gamma
    + \iO f'(y(t)) \, v\,dx
  + \iG \bigl( \fG'(\yG(t)) - \uG(t) \bigr) \, \vG\,d\Gamma ,
  \non
  \label{seconda}
\Esist
{\rm for every $v\in V$ and every $(v,\vG)\in\calV$, respectively,
and the Cauchy condition}
\Beq
  y(0) = \yz \,,\quad \yG(0) = y_{0_\Gamma } \,.
  \label{cauchy}
\Eeq
\end{definition}

\begin{remark}
{\rm It is worth noting that (recall the notation~\eqref{media})}
\Bsist
  && (\dt y(t))^\Omega = 0
  \quad \mbox{\rm for a.\,a. }\,t\in (0,T)
  \quad \mbox{\rm and }\,
  y(t)^\Omega = \mz
  \quad \mbox{{\rm for every }}\,t\in[0,T],
  \qquad
  \non
  \\
  && \quad \mbox{{\rm where }}\,\mz=(\yz)^\Omega \,\,\,
  \mbox{{\rm is the mean value of }}
  \,\yz,
  \label{conserved}
\Esist
{\rm as usual for the Cahn--Hilliard equation.}
\end{remark}

\vspace{1mm}
Now recall that $\Uad$ is a convex, closed, and bounded subset of the Banach space ${\cal X}$ and thus
contained in some bounded open ball in ${\cal X}$. For convenience, we fix such a ball once and for all,
noting that any other such ball could be used instead. The next assumption is thus rather a denotation:

\vspace{3mm}\noindent
{\bf (A4)} \quad\,The set ${\cal U}$ is some open ball in ${\cal X}$ that contains $\Uad$ and satisfies
\begin{equation}
\label{defR}
\|u_\Gamma\|_{H^1(0,T;L^2(\Gamma))}\,+\,\|u_\Gamma\|_{L^\infty(\Sigma)}\,\le\,R\quad \forall
\,u_\Gamma\in {\cal U}, 
\end{equation}
where $R>0$ is a fixed given constant.  

Concerning the well-posedness of the state sytem, we have the following result.
\Bthm
\label{daCGS}
Suppose that the general hypotheses {\bf (A1)}--{\bf (A4)} are fulfilled. Then the state system 
{\rm (\ref{Iprima})--(\ref{Icauchy})} has for any $\,u_\Gamma\in {\cal U}\,$ a unique solution $(y,\yG,w)$ in the sense of Definition \ref{DefSolu}. Moreover, there are constants $K^*_1>0$, $K_2^*>0$, and $\widetilde r_-, 
\widetilde r_+\in (r_-,r_+)$, which only depend on $\Omega$, $T$, the shape of the nonlinearities $f$ and~$\fG$,
the initial datum~$\yz$, and  the constant $R$, such that the
following holds:\\
{\rm (i)} \,\,\,Whenever $\,(y,\yG,w)$ is the solution to {\rm (\ref{Iprima})--(\ref{Icauchy})} associated with
some $u_\Gamma\in {\cal U}$ then 
\Bsist
  && \norma{(y,y_\Gamma)}_{\W{1,\infty}{\mathcal H} \cap \H1{\mathcal V} \cap \L\infty{\mathcal G}}\,+\,\norma w_{\L\infty\Hdue}\,\leq \,K_1^*\,, 
  \label{stab}
  \\[2mm]
  && \widetilde r_- \leq y \leq \widetilde r_+ \quad\,a.\,e. \,\,\,in \,\,\,Q,
  \qquad \widetilde r_- \leq y_\Gamma \leq \widetilde r_+ \quad\,a.\,e. \,\,\,on \,\,\,\Sigma.
  \label{faraway}
\Esist
{\rm (ii)} \,\,Whenever $(y_i,y_{i,\Gamma},w_i)$, $i=1,2$, are the solutions  to {\rm (\ref{Iprima})--(\ref{Icauchy})} associated with $u_{i,\Gamma}\in {\cal U}$, $i=1,2$,  then  
\begin{equation}
\label{stability}
\norma{(y_1,y_{1,\Gamma}) - (y_2,y_{2,\Gamma})}_{\H1\calH \cap \L\infty\calV}
  \,\leq\,K^*_2\, \norma{u_{1,\Gamma} - u_{2,\Gamma}}_{L^2(\Sigma)}\,.
\end{equation}
\Ethm
\Bdim
We may apply Theorems~2.2, 2.3, 2.4, 2.6 and Corollary~2.7 of \cite{CGS}
(where $\calV$ has a slightly different meaning with respect to the present paper)
to deduce that (i) holds true. Moreover, assertion (ii) is a 
consequence of \cite[Lemma 4.1]{CGSopt}.
\Edim

\vspace{1mm}
\begin{remark}
{\rm It follows from Theorem 2.3 that the control-to-state operator}
\begin{equation}
\label{ctostate} 
{\cal S}:{\cal U}\to 
W^{1,\infty}(0,T;{\cal H}) \cap H^1(0,T;{\cal V}) \cap L^\infty(0,T;{\cal G}), \quad\,
u_\Gamma\mapsto (y,y_\Gamma)\,,
  \end{equation}
{\rm is well defined and Lipschitz continuous from $\,{\cal U}$, viewed as a subset of $L^2(\Sigma)$, into
${\cal Y}$. Moreover, in view of (\ref{stab}) and (\ref{faraway}) we may assume (by possibly choosing a larger $K_1^*$) that for any $u_\Gamma\in {\cal U}$ the corresponding state
$(y,y_\Gamma)={\cal S}(u_\Gamma)$ satisfies}
\begin{equation}
\label{fbounds}
\max_{1\le i\le 4}\,\left(\|f^{(i)}(y)\|_{L^\infty(Q)}\,+\,\|f_\Gamma^{(i)}
(y_\Gamma)\|_{L^\infty(\Sigma)} \right)\,\le\,K^*_1\,.
\end{equation} 
\end{remark}

\vspace{3mm}
Next, in order to ensure the solvability of a number of linearized systems later in this paper, we introduce the
linear initial-boundary value problem
\Bsist
  & \dt \chi - \Delta \mu = 0
  \quad \hbox{in $\, Q,$}
  \label{Iprima_linear}
  \\[0.2cm]
  & \mu =  \, \dt \chi - \Delta \chi + \lambda \, \chi+g
  \quad \hbox{in $\,Q,$}
  \label{Iseconda_linear}
  \\[0.2cm]
  & \dn \mu = 0
  \quad \hbox{on $\, \Sigma,$}
  \label{Ibc_linear}
  \\[0.2cm]
  & \chi_\Gamma = \chi\suG
	\quad \hbox{on $\, \Sigma,$}
  \label{ItCorrespond_linear}
  \\[0.2cm]
  &\dt\chi_\Gamma + \dn \chi - \DeltaG\chi_\Gamma + \lambda_\Gamma \, \chi_\Gamma = g_\Gamma
  \quad \hbox{on $\, \Sigma,$}
  \label{Iterza_linear}
  \\[0.2cm]
  & \chi(0) = \chi_0
  \quad \hbox{in $\, \Omega$},\quad \chi_\Gamma(0) = \chi_{0_\Gamma}:=\chi_{0|\Gamma}
  \quad \hbox{on $\, \Gamma$}, 
  \label{Icauchy_linear}
\Esist
and its variational counterpart, namely, for almost every $t\in (0,T)$,
\begin{align}
& \iO \dt\chi(t) \, v\,dx 
  + \iO \nabla\mu(t) \cdot \nabla v\,dx = 0 \quad\mbox{for every }\,v\in V,
  \label{linprima}
  \\[1mm]
& \iO \mu(t) v\,dx
  = \iO \dt\chi(t) \, v\,dx
  + \iG \dt\chi_\Gamma(t) \, \vG\,d\Gamma
  + \iO \nabla\chi(t) \cdot \nabla v\,dx
  \non
  \\
&\quad+ \iG \nablaG\chi_\Gamma(t) \cdot \nablaG\vG\,d\Gamma
    + \iO \bigl( \lambda(t) \, \chi(t) + g(t) \bigr) \, v\,dx
  + \iG \bigl( \lambda_\Gamma(t) \, \chi_\Gamma(t) - g_\Gamma(t) \bigr) \, \vG\,d\Gamma
  \label{linseconda}\non\\[1mm]
&\qquad\mbox{for every }\,(v,\vG)\in\calV,  
\end{align}
together with the Cauchy condition
\Beq
  \chi(0) = \chi_0 ,\quad \chi_\Gamma(0) = \chi_{0_\Gamma}.
  \label{lincauchy}
\Eeq

We have the following result.
\Blem
\label{Existlin}
Suppose that $(g,g_\Gamma)\in \H1\calH\cap(L^\infty(Q)\times L^\infty(\Sigma))$ and  $(\lambda,\lambda_\Gamma)\in W^{1,\infty}(0,T;\calH)\cap (L^\infty(Q)\times L^\infty(\Sigma)) $ are given, and let $\,\chi_0\in H^2(\Omega)\,$ be 
such that $\,\chi_{0_\Gamma}
:=\chi_{0|\Gamma}\in H^2(\Gamma)$. Then the problem {\rm (\ref{Iprima_linear})--(\ref{Icauchy_linear})}
has a unique solution in the sense that there  is a unique triple  $(\chi,\chi_\Gamma,\mu)$ that fulfills
{\rm (\ref{linprima})--(\ref{lincauchy})} and whose components satisfy the analogue of the regularity requirements 
{\rm (\ref{regy})}, {\rm (\ref{regyG})}, and {\rm (\ref{regw})}, respectively.  
Moreover, there exists a constant $\,K_3^*>0$, which only depends on $\,\Omega$, $T$, $\|\lambda\|_{L^\infty(Q)}$,
and $\,\|\lambda_\Gamma\|_{L^\infty(\Sigma)}$, such that the following holds: whenever $\chi_0=0$ then
\begin{equation}
\label{eq:2.11}
\|(\chi,\chi_\Gamma)\|_{\H1\calH \cap \L\infty\calV}\,\le\,K_3^*\,\|(g,g_\Gamma)\|_{L^2(0,T;\cal H)}\,.
\end{equation}
\Elem

\Bdim
In the following, we denote by $C_i$, $i\in\mathbb N$, positive constants that only depend on the quantities mentioned in the assertion. 
First, we observe that the results 
concerning existence, uniqueness, and regularity follow from a direct application of \cite[Cor.~2.5]{CGS}. 
Now assume that $\,\chi_0=0$. Then we have $\chi^\Omega(t) = 0$ for almost every $t\in (0,T)$. 
We thus may choose in (\ref{linprima}) $\,v={\cal N}(\chi(t))$, 
and in (\ref{linseconda}) $\,v=-\chi(t)$. 
Adding the resulting equalities,
then adding two additional terms on both sides for convenience, and
integrating with respect to time, we arrive at the identity     
\begin{eqnarray}
\label{eq:2.32}
&&\frac 1 2\left(\|\chi(t)\|^2_* + \|\chi(t)\|^2_H +\| \chi_\Gamma(t)\|_{H_\Gamma}^2\right) \,+\,\int
_0^t\!\!\int_\Omega |\nabla \chi|^2\,dx\,ds\, +\,\int_0^t\!\!\int_\Gamma
 |\nabla_\Gamma \chi_\Gamma|^2 \,d\Gamma\,ds\,\nonumber\\[2mm]
&&=\int_0^t\!\!\int_\Omega \bigl( g - \lambda \, \chi \bigr) \, \chi\,dx\,ds\, +\, 
\int_0^t\!\!\int_\Gamma \bigl( g_\Gamma - \lambda_\Gamma \, \chi_\Gamma  \bigr) \, \chi_\Gamma\,d\Gamma\,ds
\nonumber
\end{eqnarray}
for all $t\in[0,T]$.
Estimating the right-hand side with the help of Young's and Poincar\'e's inequalities, 
and applying Gronwall's lemma, we have that
\begin{equation}
\label{eq:2.34}
\|(\chi,\chi_\Gamma)\|_{\L\infty\calH \cap \L2\calV} \, \le \, C_1 \, \|(g,g_\Gamma)\|_{\L2\calH}.
\end{equation}
Moreover, we may insert $\,v={\cal N}(\partial_t \chi(t))$\,
in (\ref{linprima}) and $\,v=-\partial_t \chi(t)$ in (\ref{linseconda}). 
Adding the resulting equations,
integrating with respect to time, and using (\ref{defN}), we obtain the identity
\begin{eqnarray}
\label{eq:2.35}
&&\int_0^t \|\partial_t \chi(s)\|^2_* \,ds\,+\, \int_0^t\!\!\int_\Omega |\partial_t \chi|^2\,dx\,ds\,+\,
\int_0^t\!\!\int_\Gamma |\partial_t  \chi_\Gamma|^2\,d\Gamma\,ds\,\nonumber\\[2mm]
&&+\,\frac 1 2\,(\|\nabla \chi(t)\|_H^2 \,+\,\|\nabla_\Gamma \chi_\Gamma(t)\|_{H_\Gamma}^2)\, \nonumber\\[2mm]
&&=\quad  \,\int_0^t\!\!\int_\Omega (g-\lambda \, \chi)\,\partial _t \chi \,dx\,ds\,+\,\int_0^t\!\!
\int_\Gamma (g_\Gamma-\lambda_\Gamma \, \chi_\Gamma)\,\partial _t \chi_\Gamma\,d\Gamma\,ds\,.
\end{eqnarray}
Invoking Young's inequality, we can easily infer from (\ref{eq:2.34}) and (\ref{eq:2.35}) the estimate
\begin{equation}
\label{eq:2.36}
\|(\chi,\chi_\Gamma)\|_{\H1\calH \cap \L\infty\calV} \, \le \, C_2 \, \|(g,g_\Gamma)\|_{\L2\calH},
\end{equation}
whence the assertion follows.
\Edim



\section{Differentiability properties of the control-to-state mapping}
\label{FRECHET}
\setcounter{equation}{0}
The main objective in this section is to prove that the control-to-state 
mapping is twice continuously differentiable.
 We begin our analysis with the following result.

\Bthm
\label{Fdiff}
Suppose that {\bf (A1)}--{\bf (A4)} are fulfilled. Then the following holds true:

\noindent
{\rm (i)} \,\,\,The control-to-state mapping ${\cal S}$ is \Frechet \ differentiable in ${\cal U}$
as a mapping from ${\cal U}\subset {\cal X}$ to ${\cal Y}$. 

\noindent
{\rm (ii)} \,\,For every $u_\Gamma\in {\cal U}$, the \Frechet\ derivative 
$D\calS(\uG) \in {\cal L}({\cal X},{\cal Y})$ is given as follows: for
any $\hG\in\calX$ it holds      $D\calS(\uG)\hG=(\xi,\xiG)$,
where $(\xi,\xiG,\zeta)$ with
\begin{align}
  & \xi \in \W{1,\infty}H \cap \H1V \cap \L\infty\Hdue,
  \label{reg_xi}
  \\[1mm]
  & \xiG \in \W{1,\infty}\HG \cap \H1\VG \cap \L\infty\HdueG,
  \label{reg_xiG}
  \\[1mm]
  & \zeta \in \L\infty\Hdue,
  \label{reg_zeta}
\end{align}is the unique solution to the linearized system 
\Bsist
  & \dt \xi - \Delta \zeta = 0
  \quad \hbox{in $\, Q$},
  \label{Iprima_linearized}
  \\[0.2cm]
  & \zeta =  \, \dt \xi - \Delta \xi + f''(y) \, \xi
  \quad \hbox{in $\,Q$,}
  \label{Iseconda_linearized}
  \\[0.2cm]
  & \partial_{\bf n} \zeta = 0
  \quad \hbox{on $\, \Sigma$,}
  \label{Ibc_linearized}
  \\[0.2cm]
  & \xi_\Gamma = \xi\suG
  \quad \hbox{on $\, \Sigma,$}
  \label{ItCorrespond_linearized}
  \\[0.2cm]
  &\dt\xi_\Gamma + \dn\xi_\Gamma - \DeltaG\xi_\Gamma + f_\Gamma''(\yG) \, \xi_\Gamma = h_\Gamma
  \quad \hbox{on $\, \Sigma$,}
  \label{Iterza_linearized}
  \\[0.2cm]
  & \xi(0) = 0
  \quad \hbox{in $\, \Omega$},\quad \xi_\Gamma(0) = 0
  \quad \hbox{on $\, \Gamma$} .
  \label{Icauchy_linearized}
\Esist

\noindent
{\rm (iii)} \,The mapping $D{\cal S}:{\cal U}\rightarrow{\cal L}({\cal X},{\cal Y}),\, \uG\mapsto D\calS(\uG)$, 
is Lipschitz continuous on ${\cal U}$ in the following sense: 
there is a constant $K_4^*>0$, which only depends on the data and the constant~$R$, 
such that for all $u_{1,\Gamma},u_{2,\Gamma}\in{\cal U}$ and all $h_\Gamma\in{\cal X}$ it holds
\begin{equation}
\label{stabilityLinear}
\|(D\calS(u_{1,\Gamma})-D\calS(u_{2,\Gamma}))h_\Gamma\|_{{\cal Y}}\,\leq\, 
K_4^*\,\|u_{1,\Gamma} - u_{2,\Gamma}\|_{L^2(\Sigma)} \, \|h_\Gamma\|_{L^2(\Sigma)}.
\end{equation}
\Ethm

\vspace{2mm}
\Bdim
At first, observe that the system (\ref{Iprima_linearized})--(\ref{Icauchy_linearized}) is of form (\ref{Iprima_linear})--(\ref{Icauchy_linear}), where with $(\chi,\chi_\Gamma,\mu):=(\xi,\xi_\Gamma,\zeta)$, $g\equiv 0$, $g_\Gamma:=h_\Gamma$, and $(\lambda,\lambda_\Gamma):=(f''(y),f_\Gamma''(\yG))$, the assumptions of Lemma \ref{Existlin} are fulfilled. Consequently, for every 
$h_\Gamma\in {\cal X}$, there is a unique triple $(\xi,\xi_\Gamma,\zeta)$ that 
satisfies the corresponding variational system (\ref{linprima})--(\ref{lincauchy}) and
whose components have the
regularity properties in (\ref{reg_xi}), (\ref{reg_xiG}) and (\ref{reg_zeta}).
We may therefore apply \cite[Thm. 4.2]{CGSopt} to conclude 
the validity of the assertions (i) and (ii). 

It remains to show (iii). To this end, let $\uG\in {\cal U}$ be arbitrary and let 
$k_\Gamma\in{\cal X}$ be such that $\uG+k_\Gamma\in{\cal U}$. 
We denote 
$(y^k,y^k_\Gamma)={\cal S}(\uG+k_\Gamma)$ and $(y,y_\Gamma)={\cal S}(\uG)$, 
and we assume that any $h_\Gamma\in{\cal X}$ with
$\,\|h_\Gamma\|_{\cal X}=1\,$ is given. It then suffices to show that there is some $L>0$, 
independent of $h_\Gamma$, $\uG$ and $k_\Gamma$, such that
\begin{equation}
\label{eq:3.15}
\|(\xi^k,\xi^k_\Gamma)-(\xi,\xi_\Gamma)\|_{\cal Y}\,\le\,L\,\|k_\Gamma\|_{L^2(\Sigma)}\,,
\end{equation}
where $\,(\xi^k,\xi^k_\Gamma)=D{\cal S}(\uG+k_\Gamma)h_\Gamma\,$ and 
$\,(\xi,\xi_\Gamma)=D{\cal S}(\uG)h_\Gamma$. 
For this purpose, in the following we denote by $C_i$, $i\in\mathbb N$, positive constants that neither
depend on $\uG$, $k_\Gamma$, nor on the special choice of $h_\Gamma\in{\cal X}$ with
$\,\|h_\Gamma\|_{\cal X}=1\,$. To begin with, observe that the triple 
$\,(\widehat{\xi},\widehat{\xi}_\Gamma,\widehat\zeta):=(\xi^k,\xi^k_\Gamma,\zeta^k)-(\xi,\xi_\Gamma,\zeta)$ is the unique solution to
 the variational analogue of  the initial-boundary value problem
\Bsist
  & \dt \widehat{\xi} - \Delta \widehat{\zeta} = 0
  \quad \hbox{in $\, Q,$}
  \label{Iprima_lin_Difference}
  \\[0.2cm]
  & \widehat{\zeta} =  \dt \widehat{\xi} - \Delta \widehat{\xi} + f''(y) \, 
	\widehat{\xi} + \xi^k(f''(y^k)-f''(y))
  \quad \hbox{in $\,Q$,}
  \label{Iseconda_lin_Difference}
  \\[0.2cm]
  & \dn \widehat{\zeta} = 0
  \quad \hbox{on $\, \Sigma,$}
  \label{Ibc_lin_Difference}
  \\[0.2cm]
  & \widehat{\xi}_\Gamma = \widehat{\xi}\suG\quad \hbox{on $\, \Sigma$,}
	\label{IbcCorrespond_lin_Difference}
  \\[0.2cm]
  &\dt\widehat{\xi}_\Gamma + \dn \widehat{\xi} - \DeltaG\widehat
	{\xi}_\Gamma + f_\Gamma''(y_\Gamma) \, \widehat{\xi}_\Gamma 
	= -\xi_\Gamma^k(f_\Gamma''(y_\Gamma^k)-f_\Gamma''(y_\Gamma))\quad \hbox{on $\, \Sigma$,}
  \label{Iterza_lin_Difference}
  \\[0.2cm]
  & \widehat{\xi}(0) = 0
  \quad \hbox{in $\, \Omega$},\quad \widehat{\xi}_\Gamma(0) = 0
  \quad \hbox{on $\, \Gamma$}.
  \label{Icauchy_lin_Difference}
\Esist
Moreover, the components of $\,(\widehat{\xi},\widehat{\xi}_\Gamma,\widehat\zeta)$
enjoy the regularity properties indicated in (\ref{regy}), (\ref{regyG}), and (\ref{regw}), respectively.

Now observe that it follows from Theorem 2.3, from part (i) of this proof, and from (\ref{fbounds}), that
$(g,g_\Gamma):=(\xi^k\,(f''(y^k)-f''(y)),-\xi^k_\Gamma\,(f''(y^k_\Gamma)-f''(y_\Gamma)))$
belongs to $H^1(0,T;{\cal H})\cap (L^\infty(Q)\times L^\infty(\Sigma))$, while
$(\lambda,\lambda_\Gamma):=(f''(y),f''_\Gamma(y_\Gamma))$ belongs to $W^{1,\infty}
(0,T;{\cal H})\cap (L^\infty(Q)\times L^\infty(\Sigma))$. Moreover, (\ref{fbounds}) also implies that for every $u_\Gamma\in {\cal U}$ we have for $(y,y_\Gamma)={\cal S}(u_\Gamma)$
the estimate 
\begin{equation*}
\|f''(y)\|_{L^\infty(Q)}\,+\,\|f''_\Gamma(y_\Gamma)\|_{L^\infty(\Sigma)}\,\le\,K_1^*\,.
\end{equation*}

Hence, it follows from estimate (\ref{eq:2.11}) in Lemma~\ref{Existlin} that
\begin{equation}
\label{eq:3.19}
\|(\widehat{\xi},\widehat{\xi}_\Gamma)\|_{\cal Y}\,\le\,C_1\,\left(\|\xi^k\,(f''(y^k)-f''(y))\|_{L^2(Q)} \,+\,\|\xi^k_\Gamma\,(f_\Gamma''(y^k_\Gamma)-f_\Gamma''(y_\Gamma))\|_{L^2(\Sigma)}\right)\,.
\end{equation}

\noindent Now, by the mean value theorem and (\ref{fbounds}), there exists a positive constant 
$C_2$ such that almost everywhere in $Q$ (on $\Sigma$, respectively)
\begin{equation}
\label{eq:3.20}
|f''(y^k)-f''(y)|\,\le \,C_2\,|y^k-y|\,\quad\mbox{and }\,\,\, |f_\Gamma''(y^k_\Gamma)-f_\Gamma''(y_\Gamma)|\,\le\,C_2\,|y^k_\Gamma- y_\Gamma|\,.
\end{equation}
 
\noindent At this point, we recall that ${\cal U}$ is a bounded subset of ${\cal X}$. Since
$\uG+k_\Gamma\in{\cal U}$ and $\|h_\Gamma\|_{\cal X}=1$, we thus can infer from
(\ref{fbounds}) and from
the estimate (\ref{eq:2.11}) in Lemma~\ref{Existlin} that $(\xi^k,\xi^k_\Gamma)$ is bounded in ${\cal Y}$
independently of $k_\Gamma$, $\uG$, and the choice of $h_\Gamma\in{\cal X}$ with
$\|h_\Gamma\|_{\cal X}=1$. Using the embedding $V\subset L^4(\Omega)$ and the stability
estimate proved in Theorem 2.3, we 
therefore have
\begin{eqnarray}
\label{eq:3.21}
&&\|\xi^k\,(f''(y^k)-f''(y))\|_{L^2(Q)}^2\,\le\,C_2
\int_0^T\!\!\int_\Omega\!\left(|\xi^k|^2\,|y^k- y|^2\right)\,dx\,dt\,\nonumber\\[1mm]
&&\le\,C_2\int_0^T\!\left(\|\xi^k(t)\|^2_{L^4(\Omega)}\,\|y^k(t)- y(t)\|^2_{L^4(\Omega)}\right)dt\,\nonumber \\[2mm]
&&\le\,C_3\,\|(y^k,y^k_\Gamma)- (y,y_\Gamma)\|^2_{\cal Y}
\, \le\,C_4\,\|k_\Gamma\|_{L^2(\Sigma)}^{{2}}\,.
\end{eqnarray}

\noindent Since an analogous estimate holds for the second summand in the bracket on the right-hand side of
(\ref{eq:3.19}), the assertion follows.
\Edim

\vspace{3mm}
With the Lipschitz estimate (\ref{stabilityLinear}) we are now in the position
to show the existence of the second-order \Frechet \ derivative. We have the following result.

\vspace{5mm}
\Bthm
\label{sufficient_conditions}
Assume that {\bf (A1)}--{\bf (A4)} are fulfilled. Then the following holds true: 

\noindent
{\rm (i)} \,\,\, The control-to-state operator ${\cal S}$ is twice \Frechet \ differentiable in ${\cal U}$ as a mapping from ${\cal U}\subset {\cal X}$ to ${\cal Y}$. 

\noindent
{\rm (ii)} \,\,For every $u_\Gamma\in {\cal U}$ the second Fr\'echet derivative $D^2{\cal S}(\uG)\in {\cal L}(\cal X,{\cal L}(\cal X,
\cal Y))$ is defined as follows: if $h_\Gamma, k_\Gamma\in {\cal X}$ are arbitrary then
$D^2 {\cal S}(\uG)[h_\Gamma,k_\Gamma]=:(\eta,\eta_\Gamma)$ is the unique solution to the 
initial-boundary value problem   
\Bsist
  & \dt \eta - \Delta \theta = 0
  \quad \hbox{in $\, Q,$}
  \label{Iprima_lin_second}
  \\[0.2cm]
  & \theta =  \dt \eta - \Delta \eta + f''(y) \, \eta + f^{{(3)}}(y)\, \phi \, \psi
  \quad \hbox{in $\,Q,$}
  \label{Iseconda_lin_second}
  \\[0.2cm]
  & \dn \theta = 0
  \quad \hbox{on $\, \Sigma,$}
  \label{Ibc_lin_second}
  \\[0.2cm]
  & \eta_\Gamma = \eta\suG\quad \hbox{on $\, \Sigma$,}
  \label{IbcCorrespond_lin_second}
  \\[0.2cm]
  &\dt\eta_\Gamma + \dn \eta - \DeltaG\eta_\Gamma + f_\Gamma''(\yG) \, \eta_\Gamma
  = - f_\Gamma^{{(3)}}(y_\Gamma)\,\phi_\Gamma\,\psi_\Gamma
  \quad \hbox{on $\, \Sigma$,}
  \label{Iterza_lin_second}
  \\[0.2cm]
  &\eta(0) = 0
  \quad \hbox{in $\, \Omega$},\quad \eta_\Gamma(0) = 0
  \quad \hbox{on $\, \Gamma$},
  \label{Icauchy_lin_second}
\Esist
where we have put
\begin{equation}
\label{eq:3.36}
(y,y_\Gamma)={\cal S}(\uG), \quad(\phi,\phi_\Gamma)=D{\cal S}(\uG)h_\Gamma,
\quad (\psi,\psi_\Gamma)=D{\cal S}(\uG)k_\Gamma\,.
\end{equation}

\noindent
{\rm (iii)} \,The mapping $D^2{\cal S}:{\cal U}\to {\cal L}({\cal X},{\cal L}({\cal X},{\cal Y}))$,
$\uG\mapsto D^2{\cal S}(\uG)$, is Lipschitz continuous on ${\cal U}$ in the following sense: there
exists a constant $K_5^*>0$, which only depends on the data and on the constant $R$,
 such that for every $u_{1,\Gamma},u_{2,\Gamma}\in{\cal U}$ and all $h_\Gamma, k_\Gamma\in{\cal X}$ it holds
\Beq\label{stabilitySecond}
\|(D^2\calS(u_{1,\Gamma})-D^2\calS(u_{2,\Gamma}))[h_\Gamma\,,\,k_\Gamma]\|_{{\cal Y}}\,
\leq\, K_5^*\,\|u_{1,\Gamma} - u_{2,\Gamma}\|_{L^2(\Sigma)} \, \|h_\Gamma\|_{L^2(\Sigma)} \, \|k_\Gamma\|_{L^2(\Sigma)}.
\Eeq
\Ethm

\Bdim
At first, it is easily verified that the pair $\,(g,g_\Gamma):=(f^{(3)}(y)\,\varphi\,\psi, -f^{(3)}_\Gamma(y_\Gamma)\,
\varphi_\Gamma\,\psi_\Gamma)\,$ belongs to $H^1(0,T;{\cal H})\cap (L^\infty(Q)\times L^\infty(\Sigma))$. We thus can argue as in the proof of Theorem 3.1 to deduce from Lemma 2.5 that the system (\ref{Iprima_lin_second})--(\ref{Icauchy_lin_second})
is uniquely solvable in the sense that its variational counterpart has a unique solution  
$(\eta,\eta_\Gamma,\vartheta)$ whose components enjoy the regularity indicated in
(\ref{regy}), (\ref{regyG}), and (\ref{regw}), respectively.
 Moreover, by (\ref{eq:2.11}) we have the estimate
\begin{equation}
\label{eq:3.38}
\|(\eta,\eta_\Gamma)\|_{\cal Y}\,\le\,C_1\left( \left\|f^{{(3)}}(y)\,\phi\,\psi\right\|_{L^2(Q)}\,+\,
\left\|f_\Gamma^{{(3)}}(y_\Gamma)\,\phi_\Gamma\,\psi_\Gamma\right\|_{L^2(\Sigma)}\right)\,.
\end{equation}
Here, and in the remainder of the proof of parts (i), (ii), we denote by $C_i$, $i\in\mathbb N$, positive constants that do not depend on the
quantities $h_\Gamma$, $k_\Gamma$, and $u_\Gamma$. Using (\ref{fbounds}),
and invoking the embedding $V\subset L^4(\Omega)$, we find that 
\begin{eqnarray}
\label{eq:3.39}
&&\left\|f^{{(3)}}(y)\,\phi\,\psi\right\|^2_{L^2(Q)}\,\le\,C_2
\int_0^T\!\!\int_\Omega|\phi|^2\,|\psi|^2\,dx\,dt
\,\le\,C_2\int_0^T\|\phi(t)\|^2_{L^4(\Omega)}\,\|\psi(t)\|^2_{L^4(\Omega)}\,dt\nonumber\\[2mm]
&&\le\,C_3\,\|\phi\|^2_{L^\infty(0,T;V)}\,\|\psi\|^2_
{L^\infty(0,T;V)}\,\le\,C_4\,\|h_\Gamma\|_{L^2(\Sigma)}^2
\,\|k_\Gamma\|^2_{L^2(\Sigma)}\,,
\end{eqnarray}
where the validity of the last inequality can be seen as follows: by definition 
(recall (\ref{eq:3.36})),
$(\phi,\phi_\Gamma)$ is the unique solution to the linear problem (\ref{Iprima_linearized})--(\ref{Icauchy_linearized}). We can therefore infer from (\ref{eq:2.11}) that
$\|(\phi,\phi_\Gamma)\|_{\cal Y}\,\le\,C_5\,\|h_\Gamma\|_{L^2(\Sigma)}$. By the same token, we conclude
that
$\|(\psi,\psi_\Gamma)\|_{\cal Y}\,\le\,C_6\,\|k_\Gamma\|_{L^2(\Sigma)}$. The asserted inequality
therefore follows from the definition of the norm of the space ${\cal Y}$,
and we obtain from similar reasoning that also
$$
\left\|f_\Gamma^{{(3)}}(y_\Gamma)\,\phi_\Gamma\,\psi_\Gamma\right\|_{L^2(\Sigma)}\,\le\,C_7\,\|h_\Gamma\|_{L^2(\Sigma)}\,
\|k_\Gamma\|_{L^2(\Sigma)}\,.
$$
Hence, we get
\begin{equation}
\label{eq:meanwhile}
\|(\eta,\eta_\Gamma)\|_{\cal Y}\,\le\,C_8\,\|h_\Gamma\|_{L^2(\Sigma)}\,
\|k_\Gamma\|_{L^2(\Sigma)}\,.
\end{equation}
In particular, it follows that the bilinear mapping
${\cal X}\times {\cal X}\mapsto {\cal Y}$, $[k_\Gamma,h_\Gamma]\mapsto(\eta,\eta_\Gamma)$, 
is continuous.

Now we prove the assertions concerning existence and form of the second Fr\'echet derivative. Since
${\cal U}$ is open, there is some $\Lambda>0$ such that $u_\Gamma+k_\Gamma\in {\cal U}$
whenever $\|k_\Gamma\|_{\cal X}\,\le\,\Lambda$. In the following, we only consider such perturbations $k_\Gamma\in {\cal X}$. We observe that for
$(y,y_\Gamma)={\cal S}(u_\Gamma)$ and for $(y^k,y^k_\Gamma)={\cal S}(u_\Gamma+k_\Gamma)$ the global estimates (\ref{stab})--(\ref{stability}) and (\ref{fbounds}) are
satisfied.

 After these preparations, we notice that it suffices to show that
\begin{eqnarray}
\label{eq:3.41}
&&\left\|D{\cal S}(u_\Gamma+k_\Gamma)-D{\cal S}(\uG)-D^2{\cal S}(\uG)k_\Gamma\right\|_{\cal L
({\cal X},{\cal Y})}\nonumber\\[2mm]
&&=\,\sup_{\|h_\Gamma\|_{\cal X}=1}\,
\left\|\left(D{\cal S}(u_\Gamma+k_\Gamma)-D{\cal S}(\uG)-D^2{\cal S}(\uG)k_\Gamma\right)
h_\Gamma\right\|_{\cal Y}\nonumber\\[2mm]
&&\le\,\overline C \, \norma{k_\Gamma}_{\LS2}^2
\end{eqnarray}
with a constant $\overline C$ independent of~$k_\Gamma$.

 To this end, let $h_\Gamma\in {\cal X}$ be arbitrary with $\,\|h_\Gamma\|_{{\cal X}}=1$.
We put $(\rho,\rho_\Gamma)=D{\cal S}(u_\Gamma+k_\Gamma)h_\Gamma$, define the pairs $(\phi,\phi_\Gamma),
(\psi,\psi_\Gamma)$ as in (\ref{eq:3.36}), and define
$$(\nu,\nu_\Gamma):=(\rho,\rho_\Gamma)-(\phi,\phi_\Gamma)-(\eta,\eta_\Gamma).$$
Observe that the components of  $(\nu,\nu_\Gamma)$ have the regularity properties
indicated in (\ref{regy}) and (\ref{regyG}), respectively. Moreover,
in view of (\ref{eq:3.41}), we need to show that
\begin{equation}
\label{eq:3.42}
\|(\nu,\nu_\Gamma)\|_{\cal Y}\,\le\,\overline C \, \norma{k_\Gamma}_{\LS2}^2\,.
\end{equation}
 
Now, invoking the explicit expressions for the quantities defined above, it is easily seen that the triple $(\nu,\nu_\Gamma,\pi)$ (where $\pi$ is defined below)
is the unique solution to the variational counterpart of the linear
initial-boundary value problem
\Bsist
  & \dt \nu - \Delta \pi = 0
  \quad \hbox{in $\, Q,$}
  \label{Iprima_lin_second-Differenz}
  \\[0.2cm]
  & \pi = \dt \nu - \Delta \nu + f''(y) \, \nu + \sigma 
  \quad \hbox{in $\,Q$,}
  \label{Iseconda_lin_second-Differenz}
  \\[0.2cm]
  & \dn \pi = 0
  \quad \hbox{on $\, \Sigma,$}
  \label{Ibc_lin_second-Differenz}
  \\[0.2cm]
  & \nu_\Gamma = \nu\suG
  \aand
  \dt\nu_\Gamma + \dn \nu - \DeltaG\nu_\Gamma + f_\Gamma''(\yG) \, \nu_\Gamma = \sigma_\Gamma 
  \quad \hbox{on $\, \Sigma$,}
  \label{Iterza_lin_second-Differenz}
  \\[0.2cm]
  & \nu(0) = 0
  \quad \hbox{in $\, \Omega$},\quad \nu_\Gamma(0) = 0
  \quad \hbox{on $\, \Gamma$},
  \label{Icauchy_lin_second-Differenz}
\Esist
where we have put
\begin{eqnarray}
\label{eq:3.46}
&&\sigma:=\rho\left(f''(y^k)-f''(y)\right)\,-\,f^{{(3)}}(y)\,\phi\,\psi, \nonumber\\[2mm]
&&\sigma_\Gamma:=-\rho_\Gamma\left(f_\Gamma''(y^k_\Gamma)-f_\Gamma''(y_\Gamma)\right)\,+\,f_\Gamma^{{(3)}}(y_\Gamma)
\,\phi_\Gamma\,\psi_\Gamma\,.
\end{eqnarray}
In view of (\ref{fbounds}), and since it is  easily checked that $(\sigma,\sigma_\Gamma)$ belongs to the space $H^1(0,T;{\cal H})
\linebreak\cap (L^\infty(Q))\times
L^\infty(\Sigma))$, we may again invoke the estimate (\ref{eq:2.11}) in Lemma~\ref{Existlin} 
to conclude that (\ref{eq:3.42}) is satisfied if only 
\begin{equation}
\label{eq:3.47}
\|(\sigma,\sigma_\Gamma)\|_{L^2(0,T;{\cal H})}\,\le\,\overline C \, \norma{k_\Gamma}_{\LS2}^2\,.
\end{equation} 

\noindent Applying Taylor's theorem to $f''$, and recalling (\ref{eq:3.36}), we readily see that there is a function $\,\omega_f\in L^\infty(Q)$ such that
\begin{equation}
\label{eq:3.48}
f''(y^k)-f''(y)=f^{{(3)}}(y)\,(y^k-y-\psi)\,+\,f^{{(3)}}(y)\,\psi
\,+\,\omega_f\,(y^k-y)^2\quad\mbox{a.\,e. in }\,Q\,.
\end{equation}

\noindent Hence, we have that
\begin{equation}
\label{eq:3.49}
\sigma\,=\rho\,{f^{{(3)}}}(y)\,(y^k-y-\psi)\,+\,\psi\,f^{{(3)}}(y)\,(\rho-\phi)\,+\,
\rho\,\omega_f\,(y^k-y)^2\,. 
\end{equation}

\noindent Now observe that from the proof of Fr\'echet differentiability (see inequality (4.5) in the proof of \cite[Thm. 4.2]{CGSopt}) and from 
(\ref{stabilityLinear}) we can conclude the estimates
\begin{eqnarray}
\label{eq:3.50}
&&\|(y^k,y^k_\Gamma)-(y,y_\Gamma)-(\psi,\psi_\Gamma)\|_{\cal Y}\,\le\,C_9\,\|k_\Gamma\|_{L^2(\Sigma)}^2\,,
\nonumber\\[2mm]
&&\|(\rho,\rho_\Gamma)-(\phi,\phi_\Gamma)\|_{\cal Y}\,\le\,C_{10}\,\|k_\Gamma\|_{L^2(\Sigma)}\,.
\end{eqnarray}

\noindent Moreover, we can infer from inequality (\ref{stability}) in Theorem 2.3 that
\begin{equation}
\label{eq:3.51}
\|(y^k,y^k_\Gamma)- (y,y_\Gamma)\|_{\cal Y}\,\le\,K_2^*\,\|k_\Gamma\|_{L^2(\Sigma)}\,,
\end{equation}
and it follows from Lemma~\ref{Existlin} that $(\rho,\rho_\Gamma)$ is bounded in ${\cal Y}$ by a positive 
constant that is independent of $k_\Gamma, h_\Gamma\in {\cal X}$ with
$\,\|k_\Gamma\|_{\cal X}\,\le\,\Lambda\,$ and $\,\|h_\Gamma\|_{\cal X}=1$. 

\noindent Finally, we conclude from Lemma~\ref{Existlin} (ii) that with a suitable constant $C_{11}>0$ it holds
\begin{equation}
\label{eq:3.52}
\|(\psi,\psi_\Gamma)\|_{\cal Y}\,\le\,C_{11}\,\|k_\Gamma\|_{L^2(\Sigma)}\,.
\end{equation} 

\noindent After these preparations, and invoking H\"older's inequality and the continuity of the embeddings
$V\subset L^4(\Omega)$ and $V\subset L^6(\Omega)$, we can estimate as follows:
\begin{eqnarray}
\label{eq:3.53}
&&\|\sigma\|^2_{L^2(Q)}\,\le\,C_{12}\int_0^T\!\!\int_\Omega\!\left(|\rho|^2\,|y^k-y-\psi|^2
\,+\,|\psi|^2\,|\rho-\phi|^2\,+\,|\rho|^2\,|y^k-y|^4\right)dx\,dt\,\nonumber\\[2mm]
&&{}\le\,C_{12}\!\int_0^T\!\!\!\left(\|\rho(t)\|_{L^4(\Omega)}^2\|(y^k- y-\psi)(t)\|_{L^4(\Omega)}^2
\,+\,\|\psi(t)\|_{L^4(\Omega)}^2\|\rho(t)-\phi(t)\|_{L^4(\Omega)}^2\right)dt\,\nonumber\\[2mm]
&&\quad {}+\,C_{12}\int_0^T\!\left(\|\rho(t)\|_{L^6(\Omega)}^2\,\|y^k(t)- y(t)\|^4_{L^6(\Omega)}\right)dt\,\,\nonumber\\[2mm]
&&\le\,C_{13}\,\sup_{t\in (0,T)}\left(\|\rho(t)\|_V^2\,\|(y^k- y-\psi)(t)\|_V^2\,+\,
\|\psi(t)\|_V^2\,\|\rho(t)-\phi(t)\|_V^2\right.\nonumber\\[2mm]
&&\left. \hspace*{25mm}{}+\, \|\rho(t)\|_V^2\,\|y^k(t)- y(t)\|_V^4\right)\nonumber\\[2mm]
&&\le\,C_{14}\,\|k_\Gamma\|_{L^2(\Sigma)}^4\,.
\end{eqnarray}

\noindent By the same reasoning, a similar estimate can be derived for $\|\sigma_\Gamma\|_{L^2(\Sigma)}$, which
concludes the proof of the assertions (i) and (ii).

\smallskip

Next, we prove the assertion (iii). To this end, suppose that
$\uG\in{\cal U}$ and {that $h_\Gamma$ and $k_\Gamma$ are arbitrarily chosen in $ {\cal X}$}, and let
$\delta_\Gamma\in{\cal X}$ be arbitrary with $u_\Gamma+\delta_\Gamma\in\cal U$. In the following, we will denote by $C_i$, $i\in \mathbb N$, 
positive constants that do not depend on any of these quantities.
We put
\begin{eqnarray*}
&&
(y,y_\Gamma)={\cal S}(\uG),
\quad 
(y^\delta,y^\delta_\Gamma)={\cal S}(u_\Gamma+\delta_\Gamma),
\\[1mm]
&&
(\phi,\phi_\Gamma)=D{\cal S}(\uG)h_\Gamma,
\quad 
(\phi^\delta,\phi^\delta_\Gamma)=D{\cal S}(u_\Gamma+\delta_\Gamma)h_\Gamma,
\qquad\qquad\qquad
\\[1mm]
&&
(\psi,\psi_\Gamma)=D{\cal S}(\uG)k_\Gamma,
\quad
(\psi^\delta,\psi^\delta_\Gamma)= D{\cal S}(u_\Gamma+\delta_\Gamma)k_\Gamma ,
\nonumber
\\[1mm]
&&
(\eta,\eta_\Gamma)=D^2{\cal S}(\uG)[h_\Gamma,k_\Gamma],
\quad
(\eta^\delta,\eta_\Gamma^\delta)= D^2{\cal S}(u_\Gamma+\delta_\Gamma)[h_\Gamma,k_\Gamma]\,.
\end{eqnarray*}
From the previous results, in particular, (\ref{stability}) and (\ref{stabilityLinear}),
we can infer that there is a constant $C_1>0$ such that
\begin{eqnarray}
\label{eq:3.54}
&&\|(\phi,\phi_\Gamma)\|_{\cal Y}+\|(\phi^\delta,\phi^\delta_\Gamma)\|_{\cal Y}\,\,\le\,C_1\,\|h_\Gamma\|_{L^2(\Sigma)},
\nonumber\\[2mm]
&&\|(\psi,\psi_\Gamma)\|_{\cal Y}+\|(\psi^\delta,\psi^\delta_\Gamma)\|_{\cal Y}\,\le\,C_1\,\|k_\Gamma\|_{L^2(\Sigma)},\nonumber\\[2mm]
&&\|(\eta,\eta_\Gamma)\|_{\cal Y}\,\,+\,\|(\eta^\delta,\eta^\delta_\Gamma)\|_{\cal Y}\,\,\le\,C_1\,\|h_\Gamma\|_{L^2(\Sigma)}\,\|k_\Gamma\|_{L^2(\Sigma)},\nonumber\\[2mm]
&&\|(y^\delta,y^\delta_\Gamma)-(y, y_\Gamma)\|_{\cal Y}\,\le\,C_1\,\|\delta_\Gamma\|_{L^2(\Sigma)},
\nonumber\\[2mm]
&&\|(\phi^\delta,\phi^\delta_\Gamma)-(\phi,\phi_\Gamma)\|_{\cal Y}\,\le\,C_1\,\|\delta_\Gamma\|_{L^2(\Sigma)}\, 
\|h_\Gamma\|_{L^2(\Sigma)},\nonumber\\[2mm]
&&\|(\psi^\delta,\psi^\delta_\Gamma)-(\psi,\psi_\Gamma)\|_{\cal Y}\,\le\,C_1\,
\|\delta_\Gamma\|_{L^2(\Sigma)}\,\|k_\Gamma\|_{L^2(\Sigma)}\,.
\end{eqnarray}
Now observe that $(\tilde{\eta},\tilde{\eta}_\Gamma)=(\eta^\delta,\eta^\delta_\Gamma)-(\eta,\eta_\Gamma)$ and $\tilde{\theta}=\theta^\delta-\theta$ 
(where $\vartheta^\delta$ and $\vartheta$ have their obvious meaning
corresponding to (\ref{Iseconda_lin_second})) satisfy the linear
initial-boundary value problem 
\Bsist
  & \dt \tilde{\eta} - \Delta \tilde{\theta} = 0
  \quad \hbox{in $\, Q,$}
  \label{Iprima_lin_second-Differenz-2order}
  \\[0.2cm]
  & \tilde{\theta} = \dt \tilde{\eta} - \Delta \tilde{\eta} + 
  f''(y) \, \tilde{\eta} + \sigma 
  \quad \hbox{in $\,Q$,}
  \label{Iseconda_lin_second-Differenz-2order}
  \\[0.2cm]
  & \dn \tilde{\theta} = 0
  \quad \hbox{on $\, \Sigma,$}
  \label{Ibc_lin_second-Differenz-2order}
  \\[0.2cm]
  & \tilde{\eta}_\Gamma = \tilde{\eta}\suG
  \aand
  \dt\tilde{\eta}_\Gamma + \dn \tilde{\eta} - \DeltaG\tilde{\eta}_\Gamma + f_\Gamma''(\yG) \, \tilde{\eta}_\Gamma = \sigma_\Gamma   \quad \hbox{on $\, \Sigma,$}
  \label{Iterza_lin_second-Differenz-2order}
  \\[0.2cm]
  & \tilde{\eta}(0) = 0
  \quad \hbox{in $\, \Omega$},\quad \tilde{\eta}_\Gamma(0) = 0
  \quad \hbox{on $\, \Gamma$} ,
  \label{Icauchy_lin_second-Differenz-2order}
\Esist
where we have put
\begin{eqnarray}
\label{eq:3.58}
&&\sigma=\eta^\delta(f''(y^\delta)-f''(y))+(f^{{(3)}}(y^\delta)\,\phi^\delta\,\psi^\delta
-f^{{(3)}}( y)\,\phi\,\psi)\nonumber\,,\\[1mm]
&&\sigma_\Gamma=-\eta^\delta_\Gamma
(f_\Gamma''(y^\delta_\Gamma)-f_\Gamma''( y_\Gamma))
-(f_\Gamma^{{(3)}}(y_\Gamma^\delta)\,\phi_\Gamma^\delta\,\psi_\Gamma^\delta-f_\Gamma^{{(3)}}( y_\Gamma)\,\phi_\Gamma\,\psi_\Gamma)\,.
\end{eqnarray}

\noindent The system (\ref{Iprima_lin_second-Differenz-2order})--(\ref{Icauchy_lin_second-Differenz-2order}) is again
of the form (\ref{Iprima_linear})--(\ref{Icauchy_linear}), and since 
it is readily verified that $(\sigma,\sigma_\Gamma)$ belongs to the
space $H^1(0,T;{\cal H})\cap (L^\infty(Q)\times L^\infty(\Sigma))$, we may employ  Lemma~\ref{Existlin} once more to conclude that
\begin{equation}
\label{eq:3.59}
\|(\tilde{\eta},\tilde{\eta}_\Gamma)\|_{\cal Y}\,\le\,C_2\,\|(\sigma,\sigma_\Gamma)\|_{L^2(0,T;\cal H)}\,,
\end{equation}

\noindent so that it remains to show an estimate of the form
\begin{eqnarray}
\label{eq:3.60}
&&\|(\sigma,\sigma_\Gamma)\|_{L^2(0,T;\cal H)}\,\le\,C_3\,\|\delta_\Gamma\|_{L^2(\Sigma)}\,\|h_\Gamma\|_{L^2(\Sigma)}\,\|k_\Gamma\|_{L^2(\Sigma)}\,.
\end{eqnarray}

\noindent Since
\begin{align}
\label{calc}
&f^{{(3)}}(y^\delta)\,\phi^\delta\,\psi^\delta
-f^{{(3)}}( y)\,\phi\,\psi\nonumber\\
&=\phi^\delta\,\psi\,(f^{{(3)}}(y^\delta)-f^{{(3)}}( y))+f^{{(3)}}(y^\delta)\,\phi^\delta\,(\psi^\delta-\psi)+f^{{(3)}}(y)\,\psi\,(\phi^\delta-\phi)\,,
\end{align}
\noindent we can infer from (\ref{fbounds}) that, almost everywhere in $Q$, 
\begin{equation}
\label{eq:3.61}
|\sigma|\,\le\,C_4 \,(|\eta^\delta|\,|y^\delta- y|\,+\,|\phi^\delta|\,|\psi|\,
|y^\delta- y|\,+\,|\phi^\delta|\,|\psi^\delta-\psi|\,+\,|\psi|\,|\phi^\delta-\phi|)\,.
\end{equation}

\noindent Using (\ref{eq:3.54}), H\"older's inequality, and the 
continuity of the embedding $V\subset L^4(\Omega)$, we find
\begin{eqnarray} \label{eq:3.62}
&&\int_0^T\!\!\int_\Omega\!\left(|\eta^\delta|^2\,|y^\delta-y|^2\right)
dx\,dt \,\le\,\int_0^T\! \left(\|\eta^\delta(t)\|^2_{L^4(\Omega)}
\,\|(y^\delta- y)(t)\|_{L^4(\Omega)}^2\right)dt\,\nonumber\\[2mm]
&&\le\,C_5\,\|\eta^\delta\|_{L^\infty(0,T;V)}^2\,\|y^\delta- y\|^2_
{L^\infty(0,T;V)}\,\le\,
C_6\,\|\delta_\Gamma\|_{L^2(\Sigma)}^2\,\|h_\Gamma\|_{L^2(\Sigma)}^2
\,\|k_\Gamma\|^2_{L^2(\Sigma)}\,.\nonumber\\
\end{eqnarray}

\noindent Similar reasoning yields
\begin{equation}
\label{eq:3.63}
\|\phi^\delta(\psi^\delta-\psi)\|_{L^2(Q)}^2\,+\,
\|\psi(\phi^\delta-\phi)\|_{L^2(Q)}^2\,
\le\,C_7\,\|\delta_\Gamma\|_{L^2(\Sigma)}^2\,\|h_\Gamma\|_{L^2(\Sigma)}^2
\,\|k_\Gamma\|_{L^2(\Sigma)}^2\,.
\end{equation}

\noindent Moreover, once again invoking (\ref{eq:3.54}), H\"older's inequality, and the continuity of the embedding $V\subset L^6(\Omega)$, we conclude that
\begin{eqnarray}
\label{eq:3.64}
\int_0^T\!\!\int_\Omega \!\!\left(|\phi^\delta|^2\,|\psi|^2\,|y^\delta- y|^2\right)dx\,dt\,\le
\int_0^T\!\!\!\left(\|(y^\delta- y)(t)\|_{L^6(\Omega)}^2\,\|\phi^\delta(t)\|^2_{L^6(\Omega)}\,\|\psi(t)\|^2_{L^6(\Omega)}\right)dt\,\nonumber\\[2mm]
\le\,C_8\,\|\phi^\delta\|^2_
{L^\infty(0,T;V)}\,\|\psi\|^2_{L^\infty(0,T;V)}\,\|y^\delta- y\|
^2_{L^\infty(0,T;V)}\hspace*{4.21cm}\nonumber\\[2mm]
\le\,C_9\,\|\delta_\Gamma\|_{L^2(\Sigma)}^2\,\|h_\Gamma\|_{L^2(\Sigma)}^2
\,\|k_\Gamma\|_{L^2(\Sigma)}^2\,.\hspace*{6.68cm}
\end{eqnarray}

\noindent Finally, we can estimate $\|\sigma_\Gamma\|_{L^2(\Sigma)}$, 
deriving estimates similar to (\ref{eq:3.61})--(\ref{eq:3.64}), which entails the validity of the required estimate~(\ref{eq:3.60}).
With this, the assertion is completely proved. 
\Edim


\section{Optimality conditions}
\setcounter{equation}{0}
Now that the second-order Fr\'echet-derivative of the control-to-state operator for problem {\bf (CP)} is obtained,
we can address the matter of deriving second-order sufficient optimality conditions.  
As a preparation of the corresponding theorem, we provide the adjoint system and the first-order necessary optimality conditions. 
Since these were already established in \cite{CGSopt}, we only present the results without proofs.

At first, it is easily shown (cf.\ \cite[Thm.~2.2]{CGSopt}) that {\bf (CP)}  has a solution. 
For the remainder of this paper, let us assume that $\bar u_\Gamma\in \Uad$ 
is any such minimizer and that $(\bar y,\bar y_\Gamma,\bar w)$, where
$(\bar y,\bar y_\Gamma)={\cal S}(\bar u_\Gamma)$, is the associated solution to the state system.
Recall that $(\bar y,\bar y_\Gamma,\bar w)$ has the regularity properties 
(\ref{regy}), (\ref{regyG}), and (\ref{regw}), respectively, and that 
(\ref{fbounds}) is satisfied for $(y,y_\Gamma)=(\bar y, \bar y_\Gamma)$.

The adjoint system to the problem {\bf (CP)} is formally  given by 
\Bsist
  & q + \Delta p = 0
  \quad \hbox{in $\, Q,$}
  \label{Adj1}
  \\[0.2cm]
  & -\dt (p+q) - \Delta q + f''(\bar y)\,q = \bQ(\bar y-\zQ)
  \quad \hbox{in $\,Q$,}
  \label{Adj2}
  \\[0.2cm]
  & \dn p = 0
  \quad \hbox{on $\, \Sigma,$}
  \label{Adj3}
  \\[0.2cm]
  & q_\Gamma = q\suG
  \aand
  -\dt q_\Gamma + \dn q - \Delta_\Gamma q_\Gamma + f''_\Gamma(\bar y_\Gamma)\,q_\Gamma = \bS(\bar y_\Gamma-\zS)
  \quad \hbox{on $\, \Sigma,$}
  \label{Adj4}
  \\[0.2cm]
  & (p+q)(T) = \bO(\bar y(T)-\zO)
  \quad \hbox{in $\, \Omega$},
  \label{Adj5}
	\\[0.2cm]
	& q_\Gamma(T) = \bG(\bar y_\Gamma(T)-\zG)
  \quad \hbox{on $\, \Gamma$},
  \label{Adj6}
\Esist
and was derived in \cite{CGSopt} under the additional compatibility assumption
\begin{align}\label{compatibility}
\bO=\bG=0.
\end{align}
In order to keep the technicalities at a reasonable level, we will from now on always assume that (\ref{compatibility}) is fulfilled; we remark that in \cite[Remark 5.6]{CGSopt} it has been pointed out that this assumption is dispensable at the expense of less regularity of the adjoint state variables. 

The following result was proved in \cite[Thm. 2.4]{CGSopt}.
\Bthm
Let {\bf (A1)}--{\bf (A4)} and {\rm (\ref{compatibility})} be fulfilled.
Then the adjoint system {\rm (\ref{Adj1})--(\ref{Adj6})} has a unique
solution in the following sense: there is a unique triple $(p,q,q_\Gamma)$
with the regularity properties
\Bsist
  && p \in \H1{\Hx2} \cap \L2{\Hx4},
  \label{regp}
  \\
  && q \in \H1H \cap \L2{\Hx2},
  \label{regq}
  \\
  && \qG \in \H1\HG \cap \L2{\HxG2},
  \label{regqG}
  \\
  && \qG(t) = q(t)\suG
  \quad \aat,
  \label{qqG}
	\Esist
that solves \aat\ the variational equations 
\Bsist
  && \iO q(t) \, v\,dx
  = \iO \nabla p(t) \cdot \nabla v\,dx \quad\,\forall \,v\in V,
  \label{primaadj}
  \\
  && - \iO \dt \bigl( p(t) + q(t) \bigr) \, v\,dx
  + \iO \nabla q(t) \cdot \nabla v\,dx
  + \iO f''(\yopt(t)) \, q(t) \, v \,dx
  \qquad
  \non
  \\
  && \quad {}
  - \iG \dt\qG(t) \, \vG\,d\Gamma
  + \iG \nablaG\qG(t) \cdot \nablaG\vG\,d\Gamma
  + \iG \fG''(\yGopt(t)) \, \qG(t) \, \vG\,d\Gamma
  \non
  \\
  && = \iO \bQ \bigl( \yopt(t) - \zQ(t) \bigr) v\,dx
  + \iG \bS \bigl( \yGopt(t) - \zS(t) \bigr) \vG\,d\Gamma
  \nonumber\\[2mm]
  &&\quad\mbox{for all }\,(v,\vG)\in\calV,
  \label{secondaadj}
\Esist
and the final condition
\Beq
  \iO (p+q)(T) \, v \,dx
  + \iG \qG(T) \, \vG\,d\Gamma
  = 0 \quad\,\forall (v,v_\Gamma)\in {\cal V}\,.
  \label{cauchyadj}
\Eeq
\Accorpa\Pbladj primaadj cauchyadj
\Accorpa\Regadj regp qqG
\Ethm

Now, let us introduce the 
  ``reduced cost functional'' $\redJ : \calU \to \erre$ by
\Beq
  \redJ(\uG) := \calJ(y,\yG,\uG),
  \quad \hbox{where} \quad
  (y,\yG) = \calS(\uG).
  \label{defredJ}
\Eeq
Since $\bar u_\Gamma$ is   an optimal control with associated optimal state
$(\bar y,\bar y_\Gamma)={\cal S}(\bar u_\Gamma)$, the necessary condition for optimality~is
\Beq
   D\redJ(\bar u_\Gamma)(\vG-\bar u_\Gamma)\, \geq 0,
  \quad \hbox{for every $\vG\in\Uad$},
  \label{precnopt}
\Eeq
or, written explicitly (recall that $b_\Omega=b_\Gamma=0$),
\begin{align}
\label{necessary}
&b_Q\int_0^T\!\!\int_\Omega (\bar y-z_Q)\,\xi\,dx\,dt\,+\,b_\Sigma
\int_0^T\!\!\int_\Gamma (\bar y_\Gamma-z_\Sigma)\,\xi_\Gamma\,d\Gamma\,dt
\,+\,b_0\int_0^T\!\!\int_\Gamma \bar u_\Gamma\,(v_\Gamma-\bar u_\Gamma)\,
d\Gamma\,dt \,\ge\,0\nonumber\\[2mm]
&\quad\mbox{for every }\,v_\Gamma\in\Uad, 
\end{align}
where, for any given $v_\Gamma\in\Uad$, the functions $\xi$, $\xi_\Gamma$
are the first two components of the solution triple $(\xi,\xi_\Gamma,\zeta)$
to the linearized problem (\ref{Iprima_linearized})--(\ref{Icauchy_linearized}) associated with $h_\Gamma=v_\Gamma-\bar u_\Gamma$. Moreover, since the
adjoint variables have been constructed in such a way that
\begin{equation}
\label{adjointidentity}
b_Q\int_0^T\!\!\int_\Omega (\bar y-z_Q)\,\xi\,dx\,dt\,+\,b_\Sigma
\int_0^T\!\!\int_\Gamma (\bar y_\Gamma-z_\Sigma)\,\xi_\Gamma\,d\Gamma\,dt
\,=\,\int_0^T\!\!\int_\Gamma q_\Gamma\,(v_\Gamma-\bar u_\Gamma)\,d\Gamma\,dt,
\end{equation}
we can rewrite (\ref{necessary}) in the form (see also \cite[Thm. 2.5]{CGSopt})
\Beq
  \int_0^T\!\!\int_\Gamma (\qG + \bz\, \bar u_\Gamma) (\vG - \bar u_\Gamma)\,d\Gamma\,dt \,\geq\, 0
  \quad \hbox{for every $\vG\in\Uad$}.
  \label{cnoptadj}
\Eeq
In particular, if $\bz>0$, $\bar u_\Gamma$ is the orthogonal 
projection of $-\qG/\bz$ onto $\Uad$
with respect to the standard scalar product in~$\LS2$.

After these preparations, we now derive sufficient conditions for
optimality. But, since the control-to-state operator ${\cal S}$ is not Fr\'echet differentiable on ${L^2(\Sigma)}$ 
but only on ${\cal U}\subset {\cal X}$, 
we are faced with the so-called ``two-norm discrepancy'', 
which makes it impossible to establish second-order sufficient optimality conditions 
by means of the same simple arguments as in the finite-dimensional case or, e.\,g.,
in the proof of \cite[Thm.~4.23, p.~231]{Tr}. 
It will thus be necessary to tailor the conditions in such a way as to
overcome the two-norm discrepancy. At the same time, for practical purposes the conditions should not be overly restrictive. 
For such an approach, we follow the lines of Chapter 5 in \cite{Tr}, here. Since many of the 
arguments developed here are rather similar to those employed in \cite{Tr}, we can afford to be sketchy 
and refer the reader to \cite{Tr} for
full details.

To begin with, the quadratic cost functional $\cal J$, viewed as a 
map from
 $C^0([0,T];{\cal H})\times{\cal U}$ into $\erre$, is obviously
twice continuously Fr\'echet differentiable on 
 $C^0([0,T];{\cal H})\times{\cal U}$ and thus, in particular, 
at $((\bar y,\bar y_\Gamma),\bar u_\Gamma)$. Moreover, 
 since $b_\Omega=b_\Gamma=0$, 
 we have for any 
$(( y, y_\Gamma),\uG)\in  C^0([0,T];{\cal H})\times{\cal U}$ and any $((v,v_\Gamma),h_\Gamma),
((\omega,\omega_\Gamma),k_\Gamma)\in  C^0([0,T];{\cal H})\times {\cal X}$ that
\begin{eqnarray}
\label{eq:3.65}
&&D^2{\cal J}(( y, y_\Gamma),\uG)[((v,v_\Gamma),h_\Gamma),
((\omega,\omega_\Gamma),k_\Gamma)]\nonumber\\[2mm]
&&=b_Q\int_0^T\!\!\int_\Omega v\,\omega\,dx\,dt\,+\,b_\Sigma
\int_0^T\!\!\int_\Gamma v_\Gamma\,\omega_\Gamma\,d\Gamma\,dt\,+\,
b_0\int_0^T\!\!\int_\Gamma h_\Gamma\,k_\Gamma\,d\Gamma\,dt\,.
\end{eqnarray}

\noindent It then follows from Theorem 3.2 and from the chain rule that the reduced cost functional $\tilde{{\cal J}}$
is also twice continuously Fr\'echet differentiable on ${\cal U}$. Now let  $h_\Gamma,k_\Gamma\in{\cal X}$ be arbitrary. In accordance with our previous notation,
we put
\begin{eqnarray*}
&&(\phi,\phi_\Gamma)=D{\cal S}(\bar u_\Gamma)h_\Gamma,
\quad (\psi,\psi_\Gamma)=D{\cal S}(\bar u_\Gamma)k_\Gamma,\quad
(\eta,\eta_\Gamma)=D^2{\cal S}(\bar u_\Gamma) [h_\Gamma,k_\Gamma]\,.
\end{eqnarray*}

\noindent Then a straightforward calculation resembling that carried out on page 241 in \cite{Tr}, using the 
chain rule as main tool, yields the equality
\begin{eqnarray}
\label{eq:3.66}
&&D^2\tilde{\cal J}(\bar u_\Gamma)[h_\Gamma,k_\Gamma]\,=\,D_{(y,y_\Gamma)}
{\cal J}((\bar y,\bar y_\Gamma),\bar u_\Gamma)(\eta,\eta_\Gamma)\nonumber\\[1mm]
&&\qquad+\,D^2{\cal J}((\bar y,\bar y_\Gamma),\bar u_\Gamma)
[((\phi,\phi_\Gamma),h_\Gamma)\,,\,((\psi,\psi_\Gamma),k_\Gamma)]\,.
\end{eqnarray} 
For the first summand on the right-hand side of (\ref{eq:3.66}) we have
\begin{equation}
\label{eq:3.67}
D_{(y,y_\Gamma)}{\cal J}((\bar y,\bar y_\Gamma),\bar u_\Gamma)(\eta,\eta_\Gamma)\,=\, b_Q\int_0^T\!\!\int_\Omega (\bar y-z_Q)\,\eta\,dx\,dt\,+\,b_\Sigma
\int_0^T\!\!\int_\Gamma(\bar y_\Gamma-z_\Sigma)\,\eta_\Gamma\,d\Gamma\,dt\,,
\end{equation} 
where $(\eta,\eta_\Gamma)$ solves the system (\ref{Iprima_lin_second})--(\ref{Icauchy_lin_second}). 
We now claim that
\begin{align}
\label{juerg1}
&b_Q\int_0^T\!\!\int_\Omega (\bar y-z_Q)\,\eta\,dx\,dt\,+\,b_\Sigma
\int_0^T\!\!\int_\Gamma(\bar y_\Gamma-z_\Sigma)\,\eta_\Gamma\,d\Gamma\,dt\non\\
&=-\int_0^T\!\!\int_\Omega f^{(3)}(\bar y)\,\varphi\,\psi\,q\,dx\,dt\,-\,\int_0^T\!\!\int_\Gamma
f_\Gamma^{(3)}(\bar y_\Gamma)\,\varphi_\Gamma\,\psi_\Gamma\,q_\Gamma\,d\Gamma\,dt\,.
\end{align}
To prove this claim, we test (\ref{Iprima_lin_second}) by $p$, insert $v=\theta$ in (\ref{primaadj}),
and add the resulting equations to obtain 
\begin{equation}
\label{juerg2}
0\,=\,\int_0^T\!\!\int_\Omega \left(\partial_t\eta \,p\,+\,q\,\theta\right)dx\,dt\,.
\end{equation}
Next, we test (\ref{Iseconda_lin_second}) by $q$. Since $q_{|\Gamma}=q_\Gamma$, we find the identity
\begin{align}
\label{juerg3}
&\int_0^T\!\!\int_\Omega q\,\theta\,dx\,dt \,=\,\int_0^T\!\!\int_\Omega \partial_t\eta\,q\,dx\,dt
\,+\,\int_0^T\!\!\int_\Omega \nabla\eta\cdot\nabla q\,dx\,dt 
\,+\int_0^T\!\!\int_\Gamma \partial_t\eta_\Gamma\,q_\Gamma\,d\Gamma\,dt \non\\[1mm]
&+\int_0^T\!\!\int_\Gamma \nabla_\Gamma\eta_\Gamma
\cdot \nabla_\Gamma q_\Gamma\,d\Gamma\,dt\,+
\int_0^T\!\!\int_\Omega f''(\bar y)\,\eta\,q\,dx\,dt\,+ \int_0^T\!\!\int_\Omega f^{(3)}(\bar y)\,
\varphi\,\psi\,q\,dx\,dt\non\\[1mm]
&+\int_0^T\!\!\int_\Gamma f_\Gamma''(\bar y_\Gamma)\,\eta_\Gamma\,q_\Gamma\,d\Gamma\,dt\,+
\int_0^T\!\!\int_\Gamma f_\Gamma^{(3)}(\bar y_\Gamma)\,\varphi_\Gamma\,\psi_\Gamma\,q_\Gamma\,d\Gamma\,dt\,.
\end{align} 
Now observe that the initial condition $\eta(0)=\eta_\Gamma(0)=0$ and the final condition
(\ref{cauchyadj}) imply, using integration by parts with respect to time, that 
\begin{align*}
&\int_0^T\!\!\int_\Omega \partial_t\eta \,(p+q)\,dx\,dt\,+\int_0^T\!\!\int_\Gamma \partial_t\eta_\Gamma\,q_\Gamma
\,d\Gamma\,dt\\[1mm]
&=-\int_0^T\!\!\int_\Omega  \partial_t(p+q)\,\eta\,dx\,dt\,-\int_0^T\!\!\int_\Gamma \eta_\Gamma\,\partial_t
q_\Gamma\,d\Gamma\,dt\,. 
\end{align*}
Hence, by adding (\ref{juerg2}) and (\ref{juerg3}) to each other, we obtain the identity
\begin{align}
\label{juerg4}
0&=-\int_0^T\!\!\int_\Omega \partial_t(p+q)\,\eta\,dx\,dt\,-\int_0^T\!\!\int_\Gamma \eta_\Gamma\,\partial_t q_\Gamma
\,d\Gamma\,dt\,+\int_0^T\!\!\int_\Omega \nabla\eta\cdot\nabla q\,dx\,dt\non\\[1mm]
&+\int_0^T\!\!\int_\Gamma \nabla_\Gamma\eta_\Gamma
\cdot \nabla_\Gamma q_\Gamma\,d\Gamma\,dt\,+
\int_0^T\!\!\int_\Omega f''(\bar y)\,\eta\,q\,dx\,dt\,+ \int_0^T\!\!\int_\Omega f^{(3)}(\bar y)\,
\varphi\,\psi\,q\,dx\,dt\non\\[1mm]
&+\int_0^T\!\!\int_\Gamma f_\Gamma''(\bar y_\Gamma)\,\eta_\Gamma\,q_\Gamma\,d\Gamma\,dt\,+
\int_0^T\!\!\int_\Gamma f_\Gamma^{(3)}(\bar y_\Gamma)\,\varphi_\Gamma\,\psi_\Gamma\,q_\Gamma\,d\Gamma\,dt\,.
\end{align}
Inserting $(v,v_\Gamma)=(\eta,\eta_\Gamma)$ in (\ref{secondaadj}), we finally obtain
\begin{align*}
0\,=&\int_0^T\!\!\int_\Omega\!\left(b_Q(\bar y-z_Q)\,\eta\,+\,f^{(3)}(\bar y)\,\varphi\,\psi\,q\right)dx\,dt\\[1mm]
&+\int_0^T\!\!\int_\Gamma\!\left(b_\Sigma(\bar y_\Gamma-z_\Sigma)\,\eta_\Gamma\,+\,f^{(3)}_\Gamma(\bar y_\Gamma)
\,\varphi_\Gamma\,\psi_\Gamma\,q_\Gamma\right)d\Gamma\,dt,
\end{align*}
by comparison. From this the claim (\ref{juerg1}) follows.

Now we can recall (\ref{eq:3.65})--(\ref{juerg1}) in order to find
the representation formula
\begin{align}
\label{eq:3.69}  
D^2\tilde{\cal J}(\bar u_\Gamma)[h_\Gamma,h_\Gamma]\,=\,&b_0\,\|h_\Gamma\|_{L^2(\Sigma)}^2\,+
\int_0^T\!\!\int_\Omega\!\left (b_Q-q\,f^{{(3)}}( \bar y)\right)|\phi|^2\,dx\,dt\,
\nonumber\\[1mm]
&+\int_0^T\!\!\int_\Gamma\!\left(b_\Sigma-q_\Gamma\,f_\Gamma^{{(3)}}(\bar y_\Gamma)\right)|\phi_\Gamma|^2\,d\Gamma\,dt\,
\,.
\end{align}

Equality (\ref{eq:3.69}) gives rise to hope that, under appropriate conditions, $\,D^2\tilde{\cal J}(\uopt)\,$ 
might be a positive definite operator on a suitable subset of the space $L^2(\Sigma)$.  
To formulate such a condition, we introduce
for fixed $\tau>0$ the {\em set of strongly active constraints for $\uopt$} by
\begin{eqnarray}
\label{eq:3.70}
A_\tau(\uopt)&\!\!:=\!\!&\{(x,t)\in \Sigma:\,|q_\Gamma(x,t)+\bz\,\uopt(x,t)|>\tau\}\,,
\end{eqnarray}   
and we define the 
{\em $\tau-$critical cone} $\,C_\tau(\uopt)\,$ to be the set of all $\,h_\Gamma\in{\cal X}_{M_0}:=\{h_\Gamma\in{\cal X}\,:\,\|\partial_t h_\Gamma\|_{L^2(\Sigma)}\,\leq\, M_0\}\,$
such that
\begin{eqnarray}
\label{eq:3.71}
&&h_\Gamma(x,t)\left\{
\begin{array}{lll}
=0&\mbox{if}&(x,t)\in A_\tau(\uopt)\\[1mm]
\ge 0&\mbox{if}&\uopt(x,t)=\uGmin \mbox{\, and }\,(x,t)\not \in A_\tau(\uopt)\\[1mm]
\le 0&\mbox{if}&\uopt(x,t)=\uGmax \mbox{\, and }\,(x,t)\not \in A_\tau(\uopt)
\end{array}
\right. \,.\qquad
\end{eqnarray}

\noindent After these preparations, we can formulate the second-order sufficient optimality condition (SSC) as follows:
\vspace{0.5mm}
\begin{eqnarray}
\label{eq:3.72}
&&\mbox{there exist constants }\,\delta>0\,\mbox{ and }\,\tau>0\,\mbox{ such that }\nonumber\\[2mm]
&&\qquad\qquad D^2\tilde{\cal J}(\uopt)\left[h_\Gamma,h_\Gamma\right]\,\ge\,\delta\,\|h_\Gamma\|_{L^2(\Sigma)}^2\, \quad\forall\,h_\Gamma\in C_\tau(\uopt)\,,\qquad\nonumber\\[2mm]
&&\mbox{where }\,D^2\tilde{\cal J}(\uopt)\left[h_\Gamma,h_\Gamma\right]\,
\mbox{ is given by (\ref{eq:3.69})
with }\,(\bar y,\bar y_\Gamma)={\cal S}(\uopt),\nonumber\\
&&(\phi,\phi_\Gamma)=D{\cal S}(\uopt)h_\Gamma \,\mbox{ and the associated adjoint state }\,(p,q,q_\Gamma)\,. 
\end{eqnarray}

\vspace{5mm}
The following result resembles Theorem~5.17 in \cite{Tr}.

\vspace{3mm}
\Bthm
Suppose that the conditions {\bf (A1)}--{\bf (A4)} and {\rm (\ref{compatibility})} are fulfilled,
and assume $\uopt\in\Uad$, $(\bar y ,\bar y_\Gamma)={\cal S}(\uopt)$, and that the triple
$\,(p,q,q_\Gamma)\,$ satisfies {\rm (\ref{regp})--(\ref{cauchyadj})}. Moreover, assume that the conditions 
{\rm (\ref{cnoptadj})} and {\rm (\ref{eq:3.72})} are fulfilled. Then there are constants $\varepsilon>0$ and $\sigma>0$ such that
\begin{equation}
\tilde{\cal J}(\uG)\,\ge\,\tilde{\cal J}(\uopt)\,+\,\sigma\,\|\uG - \uopt\|^2_{L^2(\Sigma)}\quad
\mbox{for all }\,\uG\in\Uad\,\mbox{ with }\,\|\uG - \uopt\|_{\cal X}\,\le\,\varepsilon\,.
\end{equation}

\noindent In particular, $\uopt$ is locally optimal for {\bf (CP)} in the sense of $\,{\cal X}$.
\Ethm
\vspace{3mm}
\Bdim
The proof closely follows that of \cite[Thm.~5.17]{Tr}, and therefore we can refer to
\cite{Tr}. We only indicate one argument that needs additional explanation. To this end, let  
$\uG\in\Uad$ be arbitrary.  Since $\tilde{\cal J}$ is twice continuously Fr\'echet 
differentiable in ${\cal U}$,
it follows from Taylor's theorem with integral remainder (see, e.\,g., \cite[Thm.~8.14.3, p.~186]{D}) that  
\begin{eqnarray}
\label{eq:3.73}
\tilde{\cal J}(\uG)-\tilde{\cal J}(\uopt)=D\tilde{\cal J}(\uopt) v_\Gamma\,+\,\frac 12\, D^2\tilde{\cal J}(\uopt)[v_\Gamma,v_\Gamma]+\,R^{\tilde{\cal J}}(\uG,\uopt)\,,
\end{eqnarray}
with $v_\Gamma=\uG-\uopt$ and the remainder 
\begin{eqnarray}
\label{eq:3.74}
R^{\tilde{\cal J}}(\uG,\uopt)=\,\int_0^1(1-s)\left(D^2\tilde{\cal J}(\uopt+s\,v_\Gamma)
-D^2\tilde{\cal J}(\uopt)\right)\left[v_\Gamma,v_\Gamma\right]\,ds.
\end{eqnarray}
Now, we estimate the integrand $(D^2\tilde{\cal J}(\uopt+s\,v_\Gamma)
-D^2\tilde{\cal J}(\uopt))[v_\Gamma,v_\Gamma]$ in (\ref{eq:3.74}). To this end, we put
\begin{eqnarray*}
&&(y^s,y^s_\Gamma)={\cal S}(\uopt+s v_\Gamma),
\quad (\phi,\phi_\Gamma)=D{\cal S}(\uopt)v_\Gamma,\quad
(\phi^s,\phi^s_\Gamma)=D{\cal S}(\uopt+s v_\Gamma)v_\Gamma,\\[1mm]
&&(\eta,\eta_\Gamma)=D^2{\cal S}(\uopt)[v_\Gamma,v_\Gamma],\quad 
(\eta^s,\eta_\Gamma^s)= D^2{\cal S}(\uopt+s v_\Gamma)
[v_\Gamma,v_\Gamma]\,,
\end{eqnarray*}
and use the representation formulas (\ref{eq:3.65})--(\ref{eq:3.67}). We obtain
\begin{align}
\label{juerg5}
D_{(y,y_\Gamma)}{\cal J}&(( y^s, y_\Gamma^s),\uopt+s v_\Gamma)(\eta^s,\eta_\Gamma^s)-D_{(y,y_\Gamma)}{\cal J}((\yopt, \yGopt),\uopt)(\eta,\eta_\Gamma)=I_1+I_{2},\qquad\nonumber \\
\end{align}
with the integrals 
\begin{align}
&I_{1}:=\,b_Q\int_0^T\!\!\int_\Omega (y^s-\yopt)\,\eta\,dx\,dt\,+\,b_\Sigma\int_0^T\!\!\int_\Gamma( y_\Gamma^s-\yGopt)\,\eta_\Gamma\,d\Gamma\,dt\,,\non\\
\label{juerg6}
&I_{2}:=\,b_Q\int_0^T\!\!\int_\Omega (y^s-z_Q)\,(\eta^s-\eta)\,dx\,dt\,+\,b_\Sigma\int_0^T\!\!\int_\Gamma
( y_\Gamma^s-z_\Sigma)\,(\eta_\Gamma^s-\eta_\Gamma)\,d\Gamma\,dt\,.
\end{align} 
Moreover, 
\begin{align}
&D^2{\cal J}(( y^s, y_\Gamma^s),\uopt+s v_\Gamma)[((\phi^s,\phi_\Gamma^s),v_\Gamma)\,,\,((\phi^s,\phi_\Gamma^s),v_\Gamma)]
\non\\[1mm]
&\quad 
-D^2{\cal J}(( \yopt, \yGopt),\uopt)[((\phi,\phi_\Gamma),v_\Gamma)\,,\,((\phi,\phi_\Gamma),v_\Gamma)]
{}=I_3,
\qquad\hbox{where}
\non\\[2mm]
\label{juerg7}
& \quad I_3:={}
b_Q\int_0^T\!\!\int_\Omega (\phi^s-\phi)(\phi^s+\phi)\,dx\,dt\,+\,b_\Sigma\int_0^T\!\!\int_\Gamma
(\phi_\Gamma^s-\phi_\Gamma)
(\phi_\Gamma^s+\phi_\Gamma)\,d\Gamma\,dt\,.
\end{align}

\noindent
We now estimate the integrals $I_{1}$, $I_2$, and $I_{3}$, where we denote by $C_i$, $i\in\enne$, constants
that neither depend on $s\in [0,1]$ nor on $\uG\in\Uad$.  At first, using the Cauchy-Schwarz inequality, we obtain
\begin{align}
\,|I_{1}|\,&\,\leq\, \max\{b_Q,b_\Sigma\}\,\|(y^s,y_\Gamma^s)-(\yopt,\yGopt)\|_{L^2(0,T;{\mathcal H})}\,\|(\eta,\eta_\Gamma)\|_{L^2(0,T;{\mathcal H})}\,\non\\[1mm]
&\,\leq\, \max\{b_Q,b_\Sigma\}\,\|(y^s,y_\Gamma^s)-(\yopt,\yGopt)\|_{\mathcal Y}\,\,\|(\eta,\eta_\Gamma)\|_{\mathcal Y}
\non\\[2mm]
\label{juerg8}
&\,\leq C_1\,s\,\|v_\Gamma\|^3_{L^2(\Sigma)},
\end{align}
where in the last inequality we have employed the estimates (\ref{stability}) and (\ref{eq:meanwhile}). Similarly,
we have
\begin{align}
\,|I_{2}|\,&\,\leq \,\max\{b_Q,b_\Sigma\}\,\|(y^s,y_\Gamma^s)-(z_Q,z_\Sigma)\|_{L^2(0,T;{\mathcal H})}\,\|(\eta^s,\eta^s_\Gamma)-(\eta,\eta_\Gamma)\|_{L^2(0,T;{\mathcal H})}\,\non\\[1mm]
&\,\leq \,\max\{b_Q,b_\Sigma\}\,\|(y^s,y_\Gamma^s)-(z_Q,z_\Sigma)\|_{L^2(0,T;{\mathcal H})}\,\,\|(\eta^s,\eta^s_\Gamma)-(\eta,\eta_\Gamma)\|_{\mathcal Y}\,\non\\[1mm]
\label{juerg9}
&\,\leq C_2\,s\,\|v_\Gamma\|^3_{L^2(\Sigma)},
\end{align}
where, for the last inequality, we used {\bf (A1)} and (\ref{stab}) to estimate the first norm and 
(\ref{stabilitySecond}) for the second one.  
Finally, we get
\begin{align}
\,|I_{3}|\,&\leq\, \max\{b_Q,b_\Sigma\}\,\|(\phi^s,\phi^s_\Gamma)-(\phi,\phi_\Gamma)\|_{L^2(0,T;{\mathcal H})}\,\|(\phi^s,\phi^s_\Gamma)+(\phi,\phi_\Gamma)\|_{L^2(0,T;{\mathcal H})}\,\non\\[1mm]
&\,\leq\, \max\{b_Q,b_\Sigma\}\,\|(\phi^s,\phi^s_\Gamma)-(\phi,\phi_\Gamma)\|_{\mathcal Y}\,\,\|(\phi^s,\phi^s_\Gamma)+(\phi,\phi_\Gamma)\|_{\mathcal Y}\,\non\\[1mm]
\label{juerg10}
&\leq C_3\,s\,\|v_\Gamma\|^3_{L^2(\Sigma)}.
\end{align}
For the last inequality we applied (\ref{stabilityLinear}) to estimate the first norm
and the triangle inequality and (\ref{eq:2.11}) to estimate the second one. 
Combining the above estimates, we thus have finally shown that 
\begin{equation}
\label{juerg11}
\left|R^{\tilde{\cal J}}(\uG,\uopt)\right|\,\,\le\,\,C_4\int_0^1(1-s)\,s\,\|v_\Gamma\|_{L^2(\Sigma)}^3\,ds\,\,\le\,\,
C_5\,\|v_\Gamma\|_{\cal X}\,\|v_\Gamma\|_{L^2(\Sigma)}^2\,,
\end{equation}
with global constants $C_4>0$ and $C_5>0$ that do not depend on the choice of $\uG\in\Uad$. 
But this means that
\begin{equation}
\frac {\left|R^{\tilde{\cal J}}(\uG,\uopt)\right|}{\|u_\Gamma-\uopt\|_{L^2(\Sigma)}^2}\,\to 0 \quad \mbox{as }\,
\|u_\Gamma-\uopt\|_{\cal X}\,\to 0.
\end{equation}
With this information at hand, we can argue along exactly the same lines as on pages 292--294 in the proof of Theorem~5.17 in \cite{Tr} to conclude the validity of the assertion.  
\Edim
\noindent \textbf{Acknowledgement.} P. Colli and G. Gilardi would like to acknowledge some financial support from the MIUR-PRIN Grant 2010A2TFX2 ``Calculus of Variations'' and the GNAMPA (Gruppo Nazionale per l'Analisi Matematica, la Probabilit\`a e le loro Applicazioni) of INdAM (Istituto Nazionale di Alta Matematica).



\vspace{3truemm}


\Begin{thebibliography}{10}

%
%
%
%
%
%
%
%
%

\bibitem{CGS}
P. Colli, G. Gilardi and J. Sprekels,
On the Cahn--Hilliard equation with dynamic boundary conditions and a dominating boundary potential,
{\em J. Math. Anal. Appl.} {\bf 419} (2014), 972-994.

\bibitem{CGSopt}
P. Colli, G. Gilardi and J. Sprekels,
A boundary control problem for the viscous Cahn--Hilliard equation with dynamic boundary conditions, 
preprint arXiv:1407.3916 [math.AP] (2014), pp. 1-27.

\bibitem{CS}
P. Colli and J. Sprekels,
Optimal control of an Allen--Cahn equation with singular potentials and dynamic boundary condition,
preprint arXiv:1212.2359~[math.AP] (2012), pp.~1-24,
to appear in SIAM J. Control Optim.

\bibitem{D}
J. Dieudonn\'e, 
Foundations of  Modern Analysis, 
Academic Press, New York, 1960.

%
%
%
%
%

\bibitem{Heik}
M. Heinkenschloss and F. Tr\"oltzsch,
Analysis of an SQP method for the control of a phase field equation,
{\it Control Cybernetics} {\bf 28} (1999) 177-211.

\bibitem{Tr}
F. Tr\"oltzsch,
{Optimal Control of Partial Differential Equations: Theory, Methods and Applications},
{\it Graduate Studies in Mathematics\/}
Vol. {\bf 112}, American Mathematical Society, Providence,
Rhode Island, 2010.

%
%
%
%

\End{thebibliography}

\End{document}

\bye